\documentclass[12pt,english,]{article}

\usepackage{bookmark}
\usepackage{graphicx}

\usepackage[round]{natbib}
\usepackage{amssymb}
\usepackage{amsfonts}
\usepackage{amsmath}
\usepackage{bbm}
\usepackage{psfrag}
\usepackage{latexsym}
\usepackage{colortbl}
\usepackage{anysize}
\usepackage{authblk}

\marginsize{2cm}{2cm}{1cm}{1cm}
\allowdisplaybreaks

\newtheorem{definition}{Definition}[section]

\newtheorem{lemma}{Lemma}[section]
\newtheorem{assumption}{Assumption}[section]
\newtheorem{assumptions}{Assumptions}[section]
\newtheorem{remark}{Remark}[section]	
\newtheorem{remarks}{Remarks}[section]
\newtheorem{example}{Example}[section]
\newtheorem{proposition}{Proposition}[section]

\providecommand{\LyX}{L\kern-.1667em\lower.25em\hbox{Y}\kern-.125emX\@}

\newcommand{\R}{\mathbb{R}}
\newcommand{\N}{\mathbb{N}}

\begin{document}

\title{GARCH-extended models: theoretical properties and applications}

\author[,a,b]{Giles-A. Nzouankeu Nana\thanks{Corresponding author, Department of Financial Mathematics, Fraunhofer Institut f\"ur Techno und Wirtschaftsmathematik (ITWM), Fraunhofer-Platz 1, 67663 Kaiserslautern, Germany, Tel.: +49-631-31600-4550, Fax: +49-631-31600-5550,  Email: giles-arnaud.nzouankeu.nana@itwm.fraunhofer.de}}
\author[a,b]{Ralf Korn}
\author[a]{Christina Erlwein-Sayer}
\affil[a]{ \footnotesize Department of Financial Mathematics, Fraunhofer Institut f\"ur Techno und Wirtschaftsmathematik.}
\affil[b]{ \footnotesize Department of Mathematics, Technische Universit\"at Kaiserslautern.}

\maketitle

\begin{abstract}
This paper is concerned with some properties of the generalized GARCH models, obtained by extending GARCH models with exogenous variables, the so-called GARCH extended (GARCHX) models. For these, we establish sufficient conditions for some properties such as stationarity, existence of moments, ergodicity, geometric ergodicity, consistence and asymptotic normality of likelihood estimators of the model parameters. For some of these properties we show that the conditions that we propose are also necessary. We further provide examples and applications to illustrate and highlight the importance of our findings.
\end{abstract}



\newpage

\section{Introduction/Motivation}

Since the development of the Autoregressive Conditional Heteroskedasticity (ARCH) model by \cite{Engl1982} and the extension to generalized ARCH (GARCH) model by \cite{Boll1986} many models of this family have been developed in oder to improve the models and to overcome some shortcomings. We can list the Exponential GARCH of \cite{Nels1991}, the GJR-GARCH of \cite{GlostJagaRunk1993}, the family GARCH (fGARCH) models of \cite{Hent1995} and the threshold GARCH of \cite{Zako1994} amongst others. The improvements proposed in all those models is done by changing the equation proposed for modeling the conditional volatility, but keeping the same variables. Since one decade another generalization ... on including exogenous variables in the model. The idea behind this procedure for financial applications is that additional sources of information help to better understand the market's behavior and hence to improve the prediction of the market's reactions. Thus one is able to earlier take dispositions against future risks. This is the case for example in \cite{BaneEtAl2011} who improve the GARCH model by introducing stock's volume as a proxy for information flow and company specific announcements in the volatility equation. \cite{SharMougKama1996} extended a GARCH(1,1) model through volume of traded stock and \cite{EnglPatt2001} introduce interest rate levels in many GARCH models. Each of these empirical studies shows a certain improvement of the GARCH model by including the exogenous process, and this improvement is presented through computational results obtained on real data. However there is only a low number of findings about the theoretical properties of such models. This is the reason of our research: we give the conditions (sufficient and sometimes also necessary) for stationarity, ergodicity, geometric ergodicity, existence of moments of the extended-GARCH, consistence and asymptotic normality of likelihood estimators. We also present some applications of these properties. As we said previously there are some pioneers in this research's direction: \cite{HanKris2012} study the asymptotic theory for the quasi maximum likelihood estimator (QMLE) in the standard GARCHX(1,1)-model with stationary and non-stationary exogenous variables. \cite{HanPark2012} analyze the asymptotic theory of QMLE of a GARCHX(1,1) with persistent covariate. Our findings generalize some of their findings, since we consider a general transformation $u$ for the exogenous variable, whereas Han and Kristensen use a special case, namely the squared function $u(x)=1+\lambda x^2$. Moreover, we not only consider the standard GARCHX model, but also almost all type of GARCHX(1,1) model. This is an important fact since many studies reveal that other GARCH type models like GJR-GARCH or E-GARCH result in a better modeling behavior than the standard GARCH model.

The rest of the paper is organized as follows: in section 2 we define the model. The third section treats the stationarity, existence of moments and ergodicity properties of our model. Further we state the conditions for the geometric ergodicity. The next section is concerned with some asymptotic theory of the likelihood estimators, namely consistence and asymptotic normality. To show that our assumptions are realistic, we provide an example which subsumes many well known GARCH models. After this we give some applications of the properties that we treat. The last section concludes the document and points out directions of future research. To improve readability all the proofs are given in the appendix.
\section{Model}

We study a general model which includes almost all GARCH(1,1) models. The model is quite similar to the one in \cite{LingMcAl2002}. However their model does not contain exogenous variables. Our model reads
\begin{eqnarray}
R_t&=&\sigma_t \varepsilon_t ,  \label{eq:first equat of the model} \\
\sigma_t^{\delta}&=& g(\varepsilon_{t-1}) + u(x_{t-1}) + c(\varepsilon_{t-1})\sigma_{t-1}^{\delta} \ , \label{eq:second equat of the model}
\end{eqnarray}

where $\{R_t\}_t$ denotes the process that we study, the exponent $\delta$ is a non-negative real value, $\sigma_t$ is the conditional volatility at time $t$, $\{\varepsilon_t\}_t$ is the noise process, $\{x_t\}_t$ represents the exogenous process used for the improvement of the modeling behavior and $g$, $c$, $u$ are real-valued non-negative continuous functions. These functions should be chosen such that $P(\sigma_t>0)=1$ for every $t$. The noise process $\{\varepsilon_t\}_t$ is supposed to be adapted to the filtration $\{\mathcal{F}_t\}_t$, where  $\mathcal{F}_t$  represents the set of all information available up to the time $t$. The relation between the noise process ${\{\varepsilon_t\}}_t$ of the GARCH and the exogenous variable $\{x_t\}_t$ is given later in the assumption \ref{assumption on noise processes}.

Our model includes many well-known GARCH (1,1) models. We can list for example (of course by setting $u \equiv 0$):
\begin{itemize}
	\item the \textit{standard GARCH(1,1)} model: $\delta=2, \ \ g(x)=\omega, \ \ c(x)=\beta+ \alpha x^2$,
	\item the \textit{GJR-GARCH(1,1) model} of \cite{GlostJagaRunk1993}: $\delta=2, \ \ g(x)=\omega, \ \ c(x)=\beta+ (\alpha + \gamma 1_{\{x > 0\}}) x^2$, where $1_{\{...\}}$ denotes the indicator function,
	\item the \textit{exponential GARCH(1,1)} of \cite{Nels1991}: $\delta \longrightarrow 0 \ \ g(x)=\omega + \alpha(\left|x\right| - E\left|x\right|), \ \ c(x)=\beta $ where $E$ denotes the expected value,
	\item the \textit{family-GARCH(1,1) (fGARCH(1,1))} models of \cite{Hent1995}: $g(x)=\omega, \ \ c(x)=\beta + \alpha (\left|x-\eta_{2j}\right|-\eta_{1j}E \left|x-\eta_{2j}\right|)$ where $\eta_{1j}, \eta_{2j}$ are non-negative real numbers,
	\item the family \textit{asymmetric power ARCH/GARCH} models of \cite{DindGranEngl1993}: $g(x)=\omega + \alpha(\left|x\right|-E\left|x\right|)^{\lambda}, \ \ c(x)=\beta$.
\end{itemize} 

We have chosen the order of the GARCH to be $(1,1)$; this is firstly due to the fact that from many studies, we can state that this order is sufficient for a good modeling (see \cite{SharMougKama1996}, \cite{BailBoll11989}, \cite{EnglPatt2001}). Secondly taking higher orders could lead to an over-fitting problem. However some results that we will present here can be easily extended to the case of higher orders. So our model can be seen as a generalization or a completion of many models.
\section{Stationarity, Ergodicity and Existence of Moments}

The properties of stationarity, ergodicity and existence of moments are three important features in the analysis of time series. Many theories and theorems only apply when the times series is stationary. The ergodicity is also an important property which implies many others and which facilitates simulations as we will see in the application presented below. Without existence of some moments it is almost impossible to get basic and also complicated results of statistics (e.g.  the law of large numbers, the central limit theorem). Hence, it is very important to find sufficient and necessary conditions under which these properties are guaranteed. We want to point out that most of the results and proofs of this part are inspired from the work of \cite{LingMcAl2002}. However their models do not contain exogenous variables. Thus our results can be considered as generalization of their findings. We also point out that we will close some gaps in the formulations and the proof of the results of \cite{LingMcAl2002}.

Let us first define the notions of stationarity and ergodicity:
\begin{definition}
\begin{enumerate}
	\item A time series $\{X_t\}_{t\in \N}$ is said to be (time-shifted) stationary if for all $n \in \N, \ t_1,...,t_n \in \N, \ h \in \N $ we have: 
$\mathcal{L}(X_{t_1},...,X_{t_n})=\mathcal{L}(X_{t_1+h},...,X_{t_n+h}).$ 
  \item A stationary process $\{X_t\}_{t\in \N}$ is called ergodic if every invariant event has either probability zero or one, where an event $A$ is called invariant if there exist a measurable set $B \in \mathfrak{B}_{\infty}$ such that for every $n \geq 1$, we have $A=\{(X_n,X_{n+1},...) \in B\}.$
\end{enumerate} 
\end{definition}

The following assumption defines the relation between the exogenous process and the noise process of the GARCH:

\begin{assumption}
The innovation process $\{\varepsilon_t\}_t$ is independent identically distributed $(i.i.d.)$ with mean $0$ and variance $1$, the exogenous process $\{x_t\}_t$ is stationary, ergodic and independent of $\{\varepsilon_t\}_t$ and the joint process $\{(\varepsilon_t, x_t)\}_t$ is adapted to $\{\mathcal{F}_t\}_t$, where $\mathcal{F}_t$ represents the set of all information available until time $t$.
\label{assumption on noise processes}
\end{assumption}

These conditions are fulfilled if one for example models the exogenous process $\{x_t\}_t$ as a measurable function of another process $\{\eta_t\}_t$ which is itself $i.i.d.$ and independent of the innovation process $\{\varepsilon_t\}_t$ or if we consider for example that $\{x_t\}_t$ is a stationary AR(1) process with a strict white noise process $\{\eta_t\}_t$ independent of $\{\varepsilon_t\}_t$, say $x_t=\varphi x_{t-1} + \eta_t$ with $\left|\varphi\right|<1$.

The first result gives us sufficient conditions under which our model is stationary and ergodic:

\begin{proposition}
Under Assumption \ref{assumption on noise processes} and $E\left|\varepsilon_t\right|^{\delta \alpha}<\infty$, $E\left(u(x_t)\right)^{\alpha}<\infty$, $E\left(c(\varepsilon_t)\right)^{\alpha}<1$, $E\left(g(\varepsilon_t)\right)^{\alpha}<\infty$ for some $\alpha \in (0,1] $, the GARCH model given by \eqref{eq:first equat of the model} and \eqref{eq:second equat of the model} admits a unique stationary solution of order $\alpha \delta$. Furthermore, the solution is strictly stationary and ergodic and admits the following causal representation: 
\begin{equation}
R_t=\sigma_t \varepsilon_t, \  \  \sigma_t^{\delta}=g(\varepsilon_{t-1})+u(x_{t-1}) + \sum\limits_{k=1}^{\infty} \left(\prod_{j=0}^{k-1} c(\varepsilon_{t-1-j})\right) \left[g(\varepsilon_{t-1-k}) + u(x_{t-1-k})\right] \ .
\label{eq: general form of the solution}
\end{equation}
\label{propoStatErgoInfty}
\end{proposition}
\begin{remark}
From this proposition \ref{propoStatErgoInfty}, we can have the surprising fact that the process $\{R_t\}_t$ is weakly stationary although the exogenous process $\{x_t\}_t$ is not. This is the case for example in the following example.
\end{remark}
\begin{example}
Let $u(x)=\left|x\right|^{1/2}$ and $\{x_t\}_t$ be iid standard Cauchy distributed. Then $\{x_t\}_t$ is not weakly stationary, since $E\left|x_t\right|=\infty$. However, proposition \ref{propoStatErgoInfty} will be valid\footnote{provided of course that the other conditions of the proposition are fulfilled.} for $\alpha \in (0,1]$. This is due to the fact that $E\left|x_t\right|^{\alpha/2}< \infty \ \ $iff$ \ \ \alpha < 2$.
\end{example}

The results stated by proposition \ref{propoStatErgoInfty} can not be used in practice since the solution requires the information until time $-\infty$, which is impossible in real life applications. In real life, we always start at some time -that we call time $0$- which is different of $-\infty$. Hence, we present below a result for the real-life process starting at some time -say 0- with a certain given distribution for $\sigma_0$. This result says that the real-life process which started at time $0$ converges with an exponential rate towards the (stationary and ergodic) process given in the above proposition and starting at time $-\infty$.

Let us denote by $\{\tilde{R}_t\}_{t\geq0}$ the same process as in equations \eqref{eq:first equat of the model} and \eqref{eq:second equat of the model}, however starting at time $0$. We denote its conditional volatility by $\tilde{\sigma}_t$, i.e.,
\begin{eqnarray}
\tilde{R}_t&=&\tilde{\sigma}_t \varepsilon_t, \\  \label{eq:first equat of the model0}
\tilde{\sigma}_t^{\delta}&=& g(\varepsilon_{t-1}) + u(x_{t-1}) + c(\varepsilon_{t-1})\tilde{\sigma}_{t-1}^{\delta}. \label{eq:second equat of the model0}
\end{eqnarray}
\begin{proposition}
Under the same assumptions of proposition \ref{propoStatErgoInfty} and given that $\tilde{\sigma}_0$ is independent of the process $\{\varepsilon_t\}_t$ and $E\left|\tilde{\sigma_0}\right|^{\delta \alpha}<\infty$, we have
\begin{equation}
\tilde{R}_t = \tilde{\sigma}_t \varepsilon_t, \  \  \  \tilde{\sigma}_t^{\delta} =\sigma_t^{\delta} + \xi_t ,
\label{eq1:propoStatErgoTime0}
\end{equation}
where $E\left|\xi_t\right|^{\alpha}=\mathcal{O}(\rho^t)$ and $0<\rho<1$ and $\sigma_t^{\delta}$ is given in proposition \ref{propoStatErgoInfty}.

Similarly,
\begin{equation}
\tilde{R}_t^{\delta}=R_t^{\delta} + \psi_t \ \ \textrm{with} \ \  E\left|\psi_t\right|^{\alpha}=\mathcal{O}(\rho^t) \ 
\label{eq2:propoStatErgoTime0}
\end{equation}
with the same $\rho$ as above and $R_t^{\delta}=\sigma_t^{\delta}\varepsilon_t^{\delta}$.
\label{propoStatErgoTime0}
\end{proposition}

This proposition states that the real-life process $\{\tilde{R}\}_t$  converges exponentially fast towards the stationary ergodic process. This fact is very helpful for practical applications: we can start our process at time $0$ at every point or with any distribution $\mathfrak{L}(\tilde{\sigma}_0)$\footnote{$\mathfrak{L}(\tilde{\sigma}_0)$ satisfying $E\left|\tilde{\sigma_0}\right|^{\delta \alpha}<\infty$ and possibly different from the stationary distribution.}, and after a small time (exponentially quickly) $t_1$ the process is nearly $\alpha \delta$-stationary and ergodic; so we can consider the process $\{\tilde{R}\}_{t \geq t_1}$ to be $\alpha \delta$-stationary and ergodic. 

Another point which is also important is to check the conditions for the existence of the moments. The following proposition gives us sufficient and necessary conditions for the existence of the $m\delta$-moment ($m\in\N$).

\begin{proposition}
Given assumption \ref{assumption on noise processes} and  $E\left|\varepsilon_t\right|^{m \delta}<\infty$, $E\left(u(x_t)\right)^{m}<\infty$, $E\left(g(\varepsilon_t)\right)^{m}<\infty$, the following assertions are equivalent:
\begin{enumerate}
	\item $E\left|R_t\right|^{m \delta} < \infty,$
	\item $E\left(c(\varepsilon_t)\right)^{m} < 1.$
\end{enumerate}
\label{existence of moments}
\end{proposition}
\section{ Geometric Ergodicity}

The geometry ergodicity is a big concern in time series. It is a very useful property from which we can derive many others, namely the mixing properties. Hence it is important to know the conditions under which our process is geometrically ergodic. For this purpose, we can adapt our situation to the one stated in \cite{Kris2005} who gives conditions for the geometric ergodicity of a class of Markov chains and who applies them to GARCH models, however without exogenous variables. So our contribution is the adaptation of their results to the context of GARCH models with exogenous variables.

We first assume that the process $\{x_t\}_t$ admits a causal finite representation, i.e. $x_t=\varphi(x_{t-1},...,x_{t-k}, \eta_t)$, where $\varphi:\mathbb{R}^{k+1}\longrightarrow \mathbb{R}$ is a continuous function and $\{\eta_t\}_t$ is an $i.i.d.\ (0,1)$ process independent of $\{\varepsilon_t\}_t$. For simplicity we will consider that k=1, since in the general case we can obtain through the representation $U_t:=(x_t,...,x_{t-k+1})=\Psi(U_{t-1}, \eta_t)$ -for a given function $\Psi$- the same result with the process $\{U_t\}_t$ instead of $\{x_t\}_t$. 
Our model \eqref{eq:first equat of the model} and \eqref{eq:second equat of the model} can be rewritten in the following way:
\begin{eqnarray*}
Y_t:=\left(\begin{array}{c} \sigma_t \\ x_{t} \end{array}\right) = \left(\begin{array}{c} \left(g(\varepsilon_{t-1}) + u(x_{t-1}) + c(\varepsilon_{t-1})\sigma_{t-1}^{\delta}\right)^{1/\delta} \\ \varphi(x_{t-1},\eta_{t}) \end{array}\right) =: G\left(Y_{t-1}, A_{t}\right),
\end{eqnarray*}
where $A_{t}=\left(\varepsilon_{t-1},\eta_{t}\right)^{'}$ is an $i.i.d$ process.

If we assume that there exists $y^{*}, a^{*} \in \mathbb{R}^2$ such that $y^{*}=G(y^{*},a^{*})$, then we can define the following quantities:
\begin{eqnarray*}
\Phi = \left.\frac{\partial G(y,a)}{\partial y}\right|_{(y^{*},a^{*})} \ ,  \  \  \ 
\Theta = \left.\frac{\partial G(y,a)}{\partial a}\right|_{(y^{*},a^{*})}. 
\end{eqnarray*}

We further define the sequence $y_n$ by $y_n=y_n(y_0,a_1,...,a_{n}) = G(y_{n-1},a_{n})$ with initial value $y_0$. We are now able to give the conditions for the geometric ergodicity of our model.

\begin{proposition}
Assume that the assumption \ref{assumption on noise processes} holds. We also assume that there exists $y^{*}$, $a^{*}$ such that $y^{*}=\lim\limits_{ n \rightarrow \infty} y_n(y_0,a^{*},...,a^{*})$. Through these values we constructed as above the matrices $\Phi$ and $\Theta$ for which we assume that $\Phi \neq 0$ and the matrix $[\Phi|\Theta]$ has full rank. We further assume that the function $G$ is differentiable at $(y^{*},a^{*})$. Finally supposing that there exists a continuous function $V: \mathbb{R}^2 \longrightarrow [1,\infty)$ satisfying $V(y)\longrightarrow \infty$ for  $\left\|y\right\|\longrightarrow \infty$ and constants $\kappa>0, \ \rho \in (0,1)$ and $T\geq1$ such that
$$E[V(Y_t)|Y_0=y] \leq \rho V(y) + \kappa \ \ \mbox{for all} \ y \in \mathbb{R}^2,$$
the process $Y_t$ is geometrically ergodic with $E[V(Y_t)] < \infty$.
\label{geometricergodicity}
\end{proposition}

Apart from the assumption of the existence of the vector $(y^{*},a^{*})$ and the function $V$, the other assumptions can be easily checked and will be certainly fulfilled by almost all classical GARCH(1,1) models. We also want to give some conditions which ensure the existence of such a function $V$. 

We assume that there exists a \textit{positive} real valued function $\bar{V}:\R^2 \longrightarrow \R,$ such that $\bar{V}(z) \longrightarrow \infty$ as $\left\|z\right\| \longrightarrow \infty$ and
\begin{equation}
\bar{V}(Y_t)\leq \alpha_1(A_{t})\bar{V}(Y_{t-1})+\beta_0(A_{t})
\end{equation}
for some function $\alpha_1$ and $\beta_0$.
 
\begin{proposition}
If such a function $\bar{V}$ exists with $E(\alpha_1(A_{t})^u)<1$, $E(\beta_0(A_t)^u)<\infty$ for a given $0\leq u <1$, then there exists a function $V$ fulfilling the desired requirements stated in the above proposition \ref{geometricergodicity}, namely $V=1+\bar{V}^{u_1}$, for some $0<u_1<u$.
\label{geoErgoExistenceOfV}
\end{proposition}


\begin{example}
We consider the T-GARCH model of \cite{Zako1994} which is a model allowing asymmetry in the volatility equation. As innovation process we consider an $i.i.d. \ \mathcal{N}(0,1)$ process and as exogenous process we set an $i.i.d.$ process. As transformation function $u$ for the exogenous variable we consider the absolute value function. The model reads

\begin{eqnarray*}
Y_t:=\left(\begin{array}{c} \sigma_t \\ x_{t} \end{array}\right) = \left(\begin{array}{c} \omega + \lambda|x_{t-1}| + (\alpha_1^{+} \varepsilon_{t-1} 1_{\{\varepsilon_{t-1}>0\}} + \alpha_1^{-} \varepsilon_{t-1} 1_{\{\varepsilon_{t-1} \leq 0\}} + \beta_1) \sigma_{t-1} \\ \gamma + \eta_{t} \end{array}\right),
\end{eqnarray*}
where $\omega,\ \alpha_1^{+},\ \alpha_1^{-},\ \beta_1>0, \  \alpha_1^{+}+\beta_1<1, \  \gamma \in \mathbb{R}_{>-1}$  and $\{\eta_t\}_t \ \ i.i.d $.

Considering the points $y^{*}=(y_1^{*},y_2^{*})^{'}=\left(\frac{\omega+\lambda(1+\gamma)}{1-(\alpha_1^{+}+\beta_1)} ,1+\gamma \right)^{'}$ and $a^{*}=\left(1,1\right)^{'}$ and the following corresponding derivatives 
\begin{equation*}
\Phi = \left.\frac{\partial G(y,a)}{\partial y}\right|_{(y^{*},a^{*})}= \left(\begin{array}{cc} \alpha_1^{+} + \beta_1 & \lambda \\ 0 & 0 \end{array}\right) \ , \  \Theta = \left.\frac{\partial G(y,a)}{\partial a}\right|_{(y^{*},a^{*})} = \left(\begin{array}{cc} \alpha_1^{+} y_1^{*} & 0 \\ 0 & 1 \end{array}\right)
\end{equation*}
at those points which fulfill all the requirements of the above proposition. This is mainly due to the fact that $\Theta$ is positive definite and $\Phi_{1,1}>0$. The only thing remaining to cover all the requirements of proposition \ref{geometricergodicity} is the existence of the function $V$. For this purpose we consider the function $\bar{V}$ as the $L^1$-norm, i.e. $\bar{V}(b_1,b_2)=|b_1|+|b_2|$. So we get\footnote{Note that $\sigma_t$ is non-negative and $c(\varepsilon_{t-1})$ is positive.}:
\begin{eqnarray*}
\bar{V}(Y_t)= \sigma_{t}+ |x_{t}| &=& \omega + \lambda|x_{t-1}| + (\alpha_1^{+} \varepsilon_{t-1} 1_{\{\varepsilon_{t-1}>0\}} + \alpha_1^{-} \varepsilon_{t-1} 1_{\{\varepsilon_{t-1}<0\}} + \beta_1) \sigma_{t-1} + |x_{t}| \\ 
&\leq& \omega + \lambda|x_{t-1}| + |x_t| + (\alpha_1^{+} \varepsilon_{t-1} 1_{\{\varepsilon_{t-1}>0\}} + \alpha_1^{-} \varepsilon_{t-1} 1_{\{\varepsilon_{t-1}<0\}} + \beta_1) (\sigma_{t-1}+ |x_{t-1}|)\\
&=& \alpha_1(A_{t})\bar{V}(Y_{t-1})+\beta_0(A_{t}),
\end{eqnarray*}
with $\alpha_1(A_{t})=c(\varepsilon_{t-1})$ and  $\beta_0(A_{t})=\omega + \lambda|x_{t-1}|+|x_t|=\omega + \lambda|\gamma + \eta_{t-1}|+ |\gamma + \eta_{t}|$.

By assuming as in the previous sections that $E(c(\varepsilon_t))^s<1$ and $E(u(x_t))^s<1$ for some $0\leq s <1$, we get the existence of $V$, namely $V=1+ \bar{V}^r$ with any $r$ satisfying $0 \leq r < s$. 

We also point out that we can apply almost the same work done in this example to almost all GARCHX models where the exogenous variable is $i.i.d.$ 
\end{example}
\section{Asymptotic Theory}

We now want to look at some asymptotic properties, namely the consistency and the asymptotic normality properties of the maximum likelihood estimator of the parameter of our model. There are some papers which treat these issues in the case of the non-extended GARCH model. Among there are \cite{LeeHans1994} who give conditions in the case of a GARCH(1,1) model and \cite{FranZako2004} who treat the general case of GARCH(p,q). There are also few works which are concerned with the case of extended standard GARCH(1,1) model like \cite{HanKris2012}, \cite{HanPark2012}. The contribution of this work is that unlike the two last papers we consider not only the standard extended GARCH(1,1) model, but a whole family of extended GARCH(1,1) model, which include almost all the GARCH-models which are usually implemented in practice. Compared to the first two papers we have included an additional variable. The proofs of this part are partly inspired by \cite{FranZako2004}, \cite{HanKris2012} and \cite{LeeHans1994}.


\subsection{Rewriting the model}
As GARCH models typically contain in their volatility equation a constant coefficient (mostly denoted by $\omega$ ) we can set: $g(x)= \omega + g_1(x)$. Additionally, we can also set $u(x)=\lambda u_1(x)$ where $\lambda$ is the coefficient of the exogenous part which has to be determined together with the rest of the coefficients. So the model reads
\begin{eqnarray}
R_t&=&\sigma_t(\theta) \varepsilon_t \ ,  \label{eq:first equat of the rewritten  model} \\
\sigma_t^{\delta}(\theta)&=&\underbrace{\omega+ g_1(\varepsilon_{t-1};\theta^{(1)})}_{g(\varepsilon_{t-1};\ \omega, \theta^{(1)})} + \underbrace{\lambda u_1(x_{t-1})}_{u(x_{t-1})} + c(\varepsilon_{t-1};\theta^{(2)})\sigma_{t-1}^{\delta}(\theta)  \ , \label{eq:second equat of the rewritten  model}
\end{eqnarray}
where $\theta^{(1)}$ (resp. $\theta^{(2)}$) represents the set of parameters contained in the function $g$ (resp. $u$), and $\theta=(\omega,\lambda,\theta^{(1)},\theta^{(2)})$ represents the set of all parameters. In order to ensure the identification of the parameters we require that $u$ is non-constant\footnote{If $u$ is a constant function, we return to the classical GARCH models without exogenous variable. Thus the requirement that $u$ is non-constant in our context is natural and trivial.}. For the construction of the likelihood function we will assume the normal distribution for the innovations process $\{\varepsilon_t\}_t$ although it can be distributed according to another distribution. This means we will construct Gaussian quasi-MLE (QMLE). However this does not mean that we are restricted to the Gaussian distribution. The model could admit any distribution, the Gaussian distribution is just used to construct the likelihood function.

The log-likelihood function (given the observations $R_1$,...,$R_n$) reads 
\begin{eqnarray}
L_n(\theta)=L_n(R_1,...,R_n;\theta)=log\left(\prod\limits_{t=1}^{n}{\frac{1}{\sqrt{2\pi\sigma_t^2}} \exp\left(-R_t^2/2\sigma_t^2\right)}\right) \ , \label{eq: unobserved log likelihood 1}
\end{eqnarray} 
where $\sigma_t$ is defined by the equation \eqref{eq: general form of the solution}. 

Ignoring the constant term\footnote{We can ignore the constant term since our goal is to maximize the likelihood function.} the likelihood function reads
\begin{equation}
L_n(\theta)=\frac{1}{n}\sum\limits_{t=1}^{n}\ell_t(\theta) \ \ \textrm{where} \ \ \ell_t(\theta)= -log(\sigma_t^2(\theta))- \frac{R_t^2}{\sigma_t^2(\theta)} \ .
\label{eq: unobserved log likelihood 2}
\end{equation}

This likelihood is however unobserved since computing the volatility $\sigma_t$ requires the knowledge of all past information until time $- \infty$ which is not the case. The observed log-likelihood function reads (ignoring the constant term)
\begin{equation}
\tilde{L}_n(\theta)= \frac{1}{n}\sum\limits_{t=1}^{n}\tilde{\ell}_t(\theta) \ \ \textrm{with} \ \ \ \tilde{\ell}_t(\theta)= -log(\tilde{\sigma}_t^2(\theta))-\frac{R_t^2}{\tilde{\sigma}_t^2(\theta)}  \ ,
\label{eq: observed log likelihood}
\end{equation}
where $\tilde{\sigma}_t$ is defined by equation \eqref{eq:second equat of the model0}\footnote{We set $\tilde{\sigma_0}^{\delta}=\omega$.}. This is the one from which the parameters should be filtered from. 

The unobserved volatility $\sigma_t$ has the big advantage that it is, under some conditions (see proposition \eqref{propoStatErgoInfty}), strictly stationary and ergodic. Thus we will use it in the proof and then show that the difference between the unobserved ($L_n(\theta)$) and the observed ($\tilde{L}_n(\theta)$) likelihood is very small.
 
We denote the parameter vector to be filtered from the model by $\theta=(\omega, \lambda,...)$. The true parameter vector is denoted by $\theta_0=(\omega_0, \lambda_0,...)$ and the parameter set will be denoted by $\Theta$. Let $\hat{\theta}_n$ be the maximum likelihood estimator given $R_1,...,R_n$, i.e.,
\begin{equation}
\hat{\theta}_n = \operatorname*{arg\,max}_{\theta \in \Theta} L_n(\theta)  \ .
\label{eq:Definition of theta}
\end{equation}   
Our aim in this chapter is to look at some asymptotic properties of $\hat{\theta}_n$.


\subsection{Consistency of the QMLE}
In this part we give conditions under which the maximum likelihood $\hat{\theta}_n$ converges\footnote{At least convergence in probability is required.} towards the true value $\theta_0$. The proofs of this part are inspired from \cite{FranZako2010} who study asymptotic properties for GARCH(p,q) models.

We assume the following:

\begin{assumptions}
\begin{enumerate}
	\item The parameter set $\Theta$ is compact, $\theta_0 \in \Theta$ and $0<\underline{\omega} \leq \omega$. The conditions of proposition \eqref{propoStatErgoInfty} are fulfilled for all elements of $\Theta$.
	\item The unobserved volatility at the true value $\theta_0$ and the innovation process $\{\varepsilon_t\}_t$ have a finite second moment i.e., $E[\sigma_t^2(\theta_0)]<\infty$ and $E[\varepsilon_t^2]<\infty$.
	\item The functions $g$ and $c$ are continuous and injective in their parameters in the sense that: 
\begin{equation*}
g(\varepsilon_t;\omega,\theta^1)=g(\varepsilon_t;\omega, \theta_0^1) \Rightarrow \omega=\omega_0 \wedge \theta^1=\theta_0^1 \textrm{ and } c(\varepsilon_t;\theta^2)=c(\varepsilon_t;\theta_0^2) \Rightarrow \theta^2=\theta_0^2. 
\end{equation*}
	\item There exists a finite set of points $\theta_i^{*}, i=1,...,n$ such that for all $\theta \in \Theta$ $g(\varepsilon_t; \theta)\leq \max\limits_{i=1,...,n}\left\{g(\varepsilon_t; \theta_i^{*})\right\}$ and $c(\varepsilon_t; \theta)\leq \max\limits_{i=1,...,n}\left\{c(\varepsilon_t; \theta_i^{*})\right\}$.
 \item $x_t|\mathfrak{F}_{t-1}$ has a non-degenerate distribution.
\end{enumerate} 
\label{ass: assumption for consistency}
\end{assumptions}

\begin{remarks} Let us mention these important points concerning the above assumptions:
\begin{itemize}
 \item We require the existence of optimum points which maximize simultaneously $u$ and $g$. However as we will see later we can have the same results if we have a finite number of optimum points for $u$ and a finite number of -possibly different- optimal points for $g$, and replace the maximum by the sum (in the assumption $(iv)$). 
  \item The fourth condition is fulfilled for almost all GARCH models. The optimal points $\theta^{*}$ are unique in the case that the functions $u$ and $g$ are linear in the parameters, but not in the general case\footnote{See for example the family fGARCH of \cite{Hent1995} in which we get four optimum points; this is due to the absolute value appearing in the volatility equation.}. 
  \item The lower bound $\underline{\omega}$ could be taken as small as possible in order to increase the set of parameters. However its positivity is required in the proofs.
  \item The third and fourth conditions seem to be complicated, but they are fulfilled by almost all usual GARCH models .
  \item Until now we have considered that the parameters of the function $u$ (in the case that $u$ has other parameters) should be fitted within the likelihood function. We can easily remark from the proofs that straightforward conditions such as imposing that the function $u$ is injective in its parameters will certainly suffice to include them in $\theta$ and to conduce the same proof.
  \item From the existence of the optimum points we can replace the first assumption $(i)$ by assuming that the conditions of proposition \ref{propoStatErgoInfty} are fulfilled for these optimum points. The same follows for all others points.
\end{itemize}
\label{remarks on consistency}
\end{remarks} 

We now state the main result of this part, namely the consistency of the quasi-maximum likelihood estimator.
\begin{proposition}
Given the Assumption \eqref{assumption on noise processes} and Assumption \eqref{ass: assumption for consistency}, a sequence of likelihood estimators $(\hat{\theta}_n)_n$ satisfying \eqref{eq:Definition of theta} converges in probability towards the true value $\theta_0$ (as $n$ goes to infinity).
\label{prop: consistency of the QMLE}
\end{proposition}

To prove this result, the following six assertions will be very helpful. The first lemma is a consequence of the assumption \ref{ass: assumption for consistency} $(iv)$ and will serve to prove among others the second lemmas. The second asserts that the difference between the observed and the unobserved volatilities is very small, that means the initial condition ($\tilde{\sigma}_0=\omega$) is asymptotically negligible. This implies the third one, which neglects the initial condition in the likelihood functions. Note that we are working under the assumptions of proposition \eqref{prop: consistency of the QMLE}

\begin{lemma} 
\begin{equation}
E(\sup\limits_{\theta\in\Theta} \sigma_t^{\delta\alpha}(\theta)) < \infty
\label{eq:first equat of lemma1 for consistence of QMLE}
\end{equation}
 and  there exist a real value $0<\rho_1<1$ such that for all $\theta \in \Theta$
\begin{equation}
E(c(\varepsilon_t; \theta)^{\alpha}) < \rho_1.
\label{eq:second equat of lemma1 for consistence of QMLE}
\end{equation}
\label{lemma: lemma1 for consistence of QMLE}
\end{lemma}
\begin{lemma}
\begin{equation}
E\left(\sup\limits_{\theta \in \Theta} \left|\sigma_t^{\delta}(\theta)- \tilde{\sigma}_t^{\delta}(\theta) \right|^{\alpha}\right) = \mathcal{O}(\rho_1^t)
\label{eq:difference between observed and unobserved volatilities}
\end{equation}
and
\begin{equation}
E\left(\sup\limits_{\theta \in \Theta} \left|\sigma_t^{2}(\theta)- \tilde{\sigma}_t^{2}(\theta) \right|^v\right) = \mathcal{O}(\rho_2^t)
\label{eq:difference between squared observed and unobserved volatilities}
\end{equation}
for some $0 < v\leq \alpha$ and $0<\rho_2<1$.
\label{lemma: difference between the both squared volatilities functions}
\end{lemma}
\begin{lemma}
\begin{equation*}
\sup\limits_{\theta \in \Theta} \left|L_n(\theta)- \tilde{L}_n(\theta) \right| \longrightarrow 0 \  \  \ \textrm{(in probability) as}  \  \ n\rightarrow \infty \ .
\end{equation*}
\label{lemma: difference between the both likelihood functions}
\end{lemma}
\begin{lemma}
\begin{equation*}
\sigma_t^{\delta}(\theta)=\sigma_t^{\delta}(\theta_0) \ \  \Rightarrow \ \ \theta=\theta_0 \ . 
\end{equation*}
\label{lemma: lemma0 for consistence}
\end{lemma}
\begin{lemma}
\begin{equation*}
E_{\theta_0}\left| \ell_t(\theta_0) \right| < \infty \textrm{ and if } \theta \neq \theta_0  \textrm{ then } E_{\theta_0}\ell_t(\theta)< E_{\theta_0}\ell_t(\theta_0) \ .
\end{equation*}
\label{lemma: lemma 2 for the proof of consistency}
\end{lemma}
\begin{lemma}
\begin{equation*}
\lim\limits_{n\rightarrow \infty} \sup\limits_{\theta^* \in \Theta} \tilde{L}_n(\theta^*) = E_{\theta_0}\ell_t(\theta_0) \  \  \ a.s. 
\end{equation*} 
\label{lemma: lemma 3 for the proof of consistency}
\end{lemma}


\subsection{Asymptotic normality of the QMLE}
Let us first recall the standard decomposition used to prove the asymptotic normality of likelihood estimators. For this purpose we denote by $$S_n(\theta)=\frac{\partial L_n}{\partial\theta} (\theta) = \frac{1}{n}\sum\limits_{t=1}^{n}\frac{\partial \ell_t(\theta)}{\partial\theta} \ \textrm{ and } \ H_n(\theta)=\frac{\partial^2 L_n}{\partial\theta\partial\theta^{'}} (\theta) = \frac{1}{n}\sum\limits_{t=1}^{n}\frac{\partial^2 \ell_t(\theta)}{\partial\theta\partial\theta^{'}} $$
the score function and the Hessian matrix of the likelihood function 
when we use the unobserved volatility $\sigma_t$; Correspondingly we denote by $\tilde{S}_n(\theta)$ and $\tilde{H}_n(\theta)$ the score function and the Hessian matrix in the case of observed volatility $\tilde{\sigma}_t$.

The derivatives of $\ell_t$ are given by\footnote{For simplicity we omit the dependence on $\varepsilon$ and $\theta$.} 
\begin{eqnarray}
\frac{\partial \ell_t}{\partial\theta_i} &=& \frac{2}{\delta} \frac{1}{\sigma_t^{\delta}} \frac{\partial \sigma_t^{\delta}}{\partial\theta_i} \left( \frac{R_t^2}{\sigma_t^{2}}-1 \right)   , \label{eq: first derivative of ell} \\
\frac{\partial^2 \ell_t}{\partial\theta_i\partial\theta_j} &=& \frac{2}{\delta} \frac{1}{\sigma_t^{\delta}} \left[-\frac{1}{\sigma_t^{\delta}} \frac{\partial \sigma_t^{\delta}}{\partial\theta_i} \frac{\partial \sigma_t^{\delta}}{\partial\theta_j} \left(\frac{2+\delta}{\delta} \frac{R_t^2}{\sigma_t^{2}}-1 \right)  +  \frac{\partial^2 \sigma_t^{\delta}}{\partial\theta_i \partial\theta_j} \left( \frac{R_t^2}{\sigma_t^{2}}-1 \right)  \right]  ,
\label{eq: second derivative of ell}
\end{eqnarray}
where the derivatives of $\sigma_t^{\delta}$ are defined as follows
\begin{eqnarray}
\frac{\partial \sigma_t^{\delta}}{\partial\theta_i} &=& \underbrace{\frac{\partial g}{\partial\theta_i}}_{M_1} + \underbrace{\frac{\partial u}{\partial\theta_i}}_{M_2} + \underbrace{\frac{\partial c}{\partial\theta_i}\sigma_{t-1}^{\delta}}_{M_3} + c\frac{\partial \sigma_{t-1}^{\delta}}{\partial\theta_i} \ ,
\label{eq: first derivative of sigma}  \\
\frac{\partial^2 \sigma_t^{\delta}}{\partial\theta_i\partial\theta_j} &=& \underbrace{\frac{\partial^2 g}{\partial\theta_i\partial\theta_j}}_{N_1} + \underbrace{\frac{\partial^2 u}{\partial\theta_i\partial\theta_j}}_{N_2} + \underbrace{\frac{\partial^2 c}{\partial\theta_i\partial\theta_j}\sigma_{t-1}^{\delta}}_{N_3} + \underbrace{\frac{\partial c}{\partial\theta_i}\frac{\partial \sigma_{t-1}^{\delta}}{\partial\theta_j}}_{N_4} + \underbrace{\frac{\partial c}{\partial\theta_j}\frac{\partial \sigma_{t-1}^{\delta}}{\partial\theta_i}}_{N_5} + c \frac{\partial^2 \sigma_{t-1}^{\delta}}{\partial\theta_i\partial\theta_j} \ .
\label{eq: second derivative of sigma} 
\end{eqnarray}

Correspondingly we get analogous expressions for the observed likelihood function $\tilde{\ell}_t$ and the observed volatility $\tilde{\sigma}_t^{\delta}$. We also set $\partial \tilde{\sigma}_0^{\delta}/\partial\theta_i=a_i$ and $\partial^2 \tilde{\sigma}_0^{\delta}/ \partial\theta_i\partial\theta_j=b_{ij}$, where $a_i,\ b_{i,j} \in \mathbb{R}$ for $i,j=1,...,n$.

Assuming that $\theta_0 \in \stackrel{\circ}{\Theta}$ and using a Taylor expansion of the score function at $\theta_0$, we get the following for the likelihood estimator $\hat{\theta}_n$ and a sufficiently large value of $n$:
\begin{eqnarray}
0=\tilde{S}_n(\hat{\theta}_n)=\tilde{S}_n(\theta_0) + \tilde{H}_n(\bar{\theta}_n)(\hat{\theta}_n-\theta_0),
\label{eq:Taylor expansion}
\end{eqnarray}
where $\bar{\theta}_n$ lies between $\hat{\theta}_n$ and $\theta_0$. We note that the first equality is due to the fact that for a sufficiently large value $n_0$, $\hat{\theta}_n, \ n \geq n_0$, is also an interior point of $\Theta$, since $\hat{\theta}_n$ converges towards $\theta_0$ and $\theta_0$ is an interior point.  

From the Taylor expansion \eqref{eq:Taylor expansion} we see that it suffices to show that 
\begin{eqnarray}
\sqrt{n}\tilde{S}_n(\theta_0) \ \longrightarrow \ \mathcal{N}(0, G(\theta_0))  \ \ \textrm{(in distribution)}, 
\label{eq: first main equation ass norm} \\
\tilde{H}_n(\bar{\theta}_n)  \ \longrightarrow \ H(\theta_0) \  \ \ \textrm{(in probability)},
\label{eq: second main equation ass norm} 
\end{eqnarray} 
with $H(\theta_0)$ non singular, in order to get by the Slutsky lemma
\begin{equation}
\sqrt{n}(\theta_n-\theta_0) \ \longrightarrow  \ \mathcal{N}(0, H^{-1}(\theta_0)G(\theta_0)(H^{-1}(\theta_0))^T);
\end{equation}
this means the asymptotic normality of the Likelihood estimators.

Let us now formulate assumptions under which we get the asymptotic normality of the QMLE.  
   
To prove the asymptotic normality of the QMLE we make the following assumptions:  
\begin{assumptions}
We assume the following:
\begin{enumerate}
	\item $\theta_0$ is an interior point of $\Theta$, i.e., $\theta_0 \in \ \stackrel{\circ}{\Theta}$.
	\item $E(\varepsilon_t^4) = \kappa < \infty$.
	\item The functions $\theta \longmapsto c(x;\theta)$ and $\theta \longmapsto g(x;\theta)$ are $ C^2(\Theta)$.
	\item There exists an $\alpha \in (0,1]$ such that for all $i,j$ and for all $\theta\in\Theta$: $$E_{\theta_0}\left( \left|\frac{\partial g(\theta)}{\partial \theta_i}\right|^{\alpha}  + \left|\frac{\partial c(\theta)}{\partial \theta_i}\right|^{\alpha} \right)<\infty \ \textrm{ and } \  E_{\theta_0}\left(\left|\frac{\partial^2 g(\theta)}{\partial \theta_i \partial \theta_j}\right|^{\alpha}  + \left|\frac{\partial^2 c(\theta)}{\partial \theta_i \partial \theta_i}\right|^{\alpha} \right)<\infty.$$
	\item For each absolute value of the derivatives of $g$ and $c$, there exists a finite number of optimum points over a neighborhood $V(\theta_0)$ of $\theta_0$. 
 \item For all $i,j,k=1,...,m$  $E_{\theta_0}\left|\partial^2\ell_t(\theta)/\partial\theta_i \partial\theta_j\right|<\infty$ and  $E_{\theta_0}\sup \limits_{\theta \in V^{'}(\theta_0)} \left|\partial^3\ell_t(\theta)/\partial\theta_i \partial\theta_j \partial\theta_k\right|<\infty$ for some neighborhood $V^{'}(\theta_0)$ of $\theta_0$.
		\item $A=E_{\theta_0}\left(\frac{1}{\sigma_t^{2\delta}}\frac{\partial\sigma_t^{\delta}(\theta_0)}{\partial \theta} \frac{\partial\sigma_t^{\delta}(\theta_0)}{\partial \theta^{'}}\right)$ exists and is invertible.
\end{enumerate}
\label{ass: assumptions for the ass norm}
\end{assumptions}

\begin{remarks}
\begin{itemize}
	\item The first assumption is natural since if $\theta$ lies on the frontier of $\Theta$ we would have at most that $\sqrt{n}(\hat{\theta}_n-\theta_0)$ is a truncated normal distributed.
  \item The second condition is due to the fact the Gaussian likelihood involves the square of the innovation process; Hence its fourth moment should exist to ensure the existence of the variance of the score function (see equation \eqref{eq: first main equation ass norm}).
  \item The differentiability condition made in Assumption 3 is imposed to ensure the existence of the first and the second derivative of the volatility process. The same holds true for assumptions 4, which are similar to those made in proposition \ref{propoStatErgoInfty}. They guarantee -as shown in proposition \ref{propoStatErgoInfty}- that the first and second derivatives of $\sigma_t^{\delta}$ are stationary and ergodic. 
  \item We note that the existence of optimum points made in assumption 4 is reasonable for many GARCH models as we will see later in the example. 
  We also note that the optimum points for one function are allowed to be different to the ones of another function.
\end{itemize}
\end{remarks}

\begin{proposition}
Under the assumptions \ref{ass: assumption for consistency}, \ref{ass: assumptions for the ass norm} and those of proposition \ref{propoStatErgoInfty}, the likelihood is asymptotically normal distributed with mean $\theta_0$, concretely,
\begin{equation}
\sqrt{n}(\hat{\theta}_n-\theta_0) \ \longrightarrow  \ \mathcal{N}\left(0,\frac{\delta^2}{4}(\kappa-1)A^{-1}\right).
\label{eq:ass norm}
\end{equation}
\label{prop: ass norm of QMLE}
\end{proposition}
To prove the above proposition \ref{prop: ass norm of QMLE}, it suffices to prove the equations \eqref{eq: first main equation ass norm} and \eqref{eq: second main equation ass norm}. For this purpose we will show the following lemmas:

\begin{lemma}
$$var\left(\frac{\partial \ell_t}{\partial \theta }(\theta_0)\right) = \frac{4}{\delta^2}(\kappa-1)A. $$
\label{lemma: first lemma for ass norm}
\end{lemma}
\begin{lemma} For all $i,j=1,...,m$ the process $\partial \ell_t/\partial\theta_i$ (resp. $\partial^2 \ell_t/\partial\theta_i\partial\theta_j$) is $\alpha/2$ (resp. $\alpha/4$) stationary and ergodic.
\label{lemma: second lemma for ass norm}
\end{lemma}
\begin{lemma}
\begin{eqnarray}
\frac{\partial\sigma_t^{\delta}(\theta)}{\partial\theta} &=&\frac{\partial\tilde{\sigma}_t^{\delta}(\theta)}{\partial\theta} + \eta_{1,t}(\theta) \  \  \ \textrm{ where }  \  \ \  E_{\theta_0}\left(\sup\limits_{\theta\in V(\theta_0)}|\eta_{1,t}(\theta)|^{\alpha/2}\right)= \mathcal{O}(\rho_3^t),  
\label{eq: diff observ and unobserv first der volatility} \\
\frac{\partial^2 \sigma_t^{\delta}(\theta)}{\partial\theta\partial\theta^{'}} &=& \frac{\partial^2 \tilde{\sigma}_t^{\delta}(\theta)} {\partial\theta\partial\theta^{'}}  + \eta_{2,t}(\theta) \  \  \  \textrm{ where }  \  \  \ E_{\theta_0}\left(\sup\limits_{\theta\in V(\theta_0)} |\eta_{2,t}(\theta)|^{\alpha/4}\right)=\mathcal{O}(\rho_4^t),
\label{eq: diff observ and unobserv second der volatility}
\end{eqnarray}
with $0 \leq \rho_3, \rho_4 < 1$.

There also exist some real numbers $0 \leq v < \alpha/2$ such that the difference of the squared volatilities -at the power $v$- are asymptotically exponentially negligible, i.e.,
\begin{equation}
\sigma_t^{2}(\theta) = \tilde{\sigma}_t^{2}(\theta) + \eta_{3,t}(\theta) \  \  \ \textrm{ where }  \  \  \  E_{\theta_0}\left(\sup\limits_{\theta\in V(\theta_0)}|\eta_{3,t}(\theta)|^{v}\right)=\mathcal{O}(\rho_5^t) \ ,
\label{eq: diff observ and unobserv square first der volatility}
\end{equation}
with $0\leq\rho_5<1$.
\label{lemma: third lemma for ass norm}
\end{lemma}
\begin{lemma}
\begin{eqnarray}
\left\|\sqrt{n}\left(S_n(\theta_0)- \tilde{S}_n(\theta_0)\right)\right\| \ \longrightarrow \ 0 \ \ (\textrm{in prob.})  \  \  as  \  \  n\rightarrow\infty,
\label{eq:diff Sn and tilde Sn} \\
\sup\limits_{\theta \in V(\theta_0)}\left\|H_n(\theta)- \tilde{H}_n(\theta)\right\| \ \longrightarrow \ 0 \ \ (\textrm{in prob.})  \  \  as  \  \  n\rightarrow\infty.
\label{eq:diff Hn and tilde Hn}
\end{eqnarray}
\label{lemma: fourth lemma for ass norm diff ell and tilde ell}
\end{lemma}

\section{Example}
We will apply our theory to study a whole family of GARCH models, namely the fGARCH family of \cite{Hent1995}, which subsumes almost all well known GARCH models (standard GARCH, absolute value GARCH, GJR GARCH, Threshold GARCH, non linear ARCH, non-linear asymmetric GARCH, asymmetric power ARCH, full fGARCH,...) and also some subfamily as we will see later. Therefore we look at
\begin{equation*}
R_t=\sigma_t \varepsilon_t, \ \ \ \sigma_t^{\delta}=\underbrace{\omega}_{g(\varepsilon_{t-1};\theta)}+ \underbrace{\lambda u(x_{t-1})}_{u(x_{t-1};\theta)}+ \underbrace{\left[\alpha_1\left(\left| \varepsilon_{t-1} - \eta_2 \right|-  \eta_1 \left( \varepsilon_{t-1} - \eta_2 \right)\right)^{\gamma} + \beta_1 \right]}_{c(\varepsilon_{t-1};\theta)} \sigma_{t-1}^{\delta} \ ,
\label{eq:fGARCH}
\end{equation*}
where the innovation process $\{\varepsilon_t\}_t$ is an $i.i.d.$ process which admits a moment of order four. It should also be non degenerate and should attain both, positive and negative values with positive probabilities. This could be for example the $\mathcal{N}(0,1)$ distribution or a \textit{Student}-$t_{\nu}$ with $\nu>4$ degrees of freedom. $\{\varepsilon_t\}_t$ is independent of $\{x_t\}_t$ which is also $i.i.d.$. The function $u$ should be non constant. It follows by supposing that $x_t$ is non degenerate that $u(x_t)$ is also non degenerate. 
\begin{eqnarray*}
\Theta=\{(\omega,\lambda, \alpha_1, \beta_1, \eta_1, \eta_2) \ | \ 0 < \underline{\omega} \leq \omega \leq \overline{\omega} < \infty, \ -1 < \underline{\eta_1} \leq \eta_1  \leq \overline{\eta_1} < 1, \ -\infty < \underline{\eta_2} \leq \eta_2 \leq \overline{\eta_2} < \infty, \\
\ 0 < \underline{\lambda_1} \leq \lambda_1 \leq \overline{\lambda_1} < \infty, \ 0 \leq \underline{\beta_1} \leq \beta_1 \leq \overline{\beta_1} < 1, \ 0 < \underline{\alpha_1} \leq \alpha_1 \leq \overline{\alpha_1} < \infty\}
\end{eqnarray*}

\subsection{Stationarity and Ergodicity (Prop. \ref{propoStatErgoInfty})}

We suppose that $E(u_1(x_t)^{\alpha})<\infty$. We set $U_{\alpha}=\left(|\varepsilon_{t-1} - \eta_2|- \eta_1( \varepsilon_{t-1} - \eta_2)\right)^{\gamma}$.  Since $\varepsilon_t$ admits a fourth moment, we have $E(\varepsilon_t)^{\delta\alpha}<\infty$ for all $\delta>0$ and for all $\alpha\in (0,1]$ such that $\delta\alpha \leq 4$. 
$E(g(\varepsilon_{t-1};\theta)^{\alpha})=E(\omega^{\alpha})= \omega^{\alpha} <\infty$ for all $\alpha\in (0,1]$.
$E(c(\varepsilon_{t-1};\theta)^{\alpha}) =  E\left(\alpha_1 U_{\alpha} + \beta_1\right)^{\alpha}$. So we can require that $E\left(\alpha_1 U_{\alpha} + \beta_1\right)^{\alpha}<1$ for a given $\alpha$ to get the $\delta \alpha$-stationarity and ergodicity of our GARCH model.

In particular we have $\beta_1<1$.

\subsection{Existence of moments (Prop. \ref{existence of moments})}

$E(g(\varepsilon_{t-1};\theta)^{m})=E(\omega^{m})= \omega^{m} <\infty$ for all $m \in \mathbb{N}$.
We suppose that $E(u_1(x_t)^{m})<\infty$ and $E(\varepsilon_t)^{m \delta}<\infty$.
So we have $$E(R_t^{m\delta})<\infty  \hspace{0.4cm}  \textrm{iff} \hspace{0.4cm}  E(\sigma_t^{m\delta}(\theta))<\infty \hspace{0.4cm}  \textrm{iff} \hspace{0.4cm}  E\left[ \alpha_1 \left(|\varepsilon_{t-1} - \eta_2|- \eta_1( \varepsilon_{t-1} - \eta_2)\right)^{\gamma} + \beta_1 \right]^m < 1 .$$

\subsection{Asymptotic Theory}

\subsubsection{Consistency (Prop. \ref{prop: consistency of the QMLE})}

We check whether the the assumptions \ref{ass: assumption for consistency} are fulfilled. 
 
$1.:$ We suppose there exist an $\alpha \in (0,1]$ such that for all $\theta \in \Theta$ we have $E\left(\alpha_1 U_{\alpha}(\theta) + \beta_1 \right)^{\alpha} <1$ where $U_{\alpha}(\theta)=\left(|\varepsilon_{t-1} - \eta_2|- \eta_1( \varepsilon_{t-1} - \eta_2)\right)^{\gamma}$. However due to the existence of optimum points listed below, we only need to check this above condition for the optimum points and not for all elements $\theta\in \Theta$. This means we require that this conditions holds only at the optimum points defined below.

$2.:$ The innovation process $\varepsilon_t$ admits moments of fourth order, thus also of second order. We suppose that the square of the volatility  at $\theta_0= (\omega^0, \lambda^0, \alpha_1^0, \beta_1^0, \eta_1^0, \eta_2^0)$ exists, this means from the above findings that it suffices to have $E(u(x_t)^m)<\infty$ and \\ $E\left[ \alpha_1^0 \left(|\varepsilon_{t-1} - \eta_2^0|- \eta_1^0 (\varepsilon_{t-1} - \eta_2^0)\right)^{\gamma} + \beta_1^0 \right]^m<1$, where $m$ is chosen such that $m\delta \geq 2$.

$3.:$: Given two points $\hat{\theta}=(\hat{\omega},\hat{\lambda},\hat{\alpha}_1, \hat{\beta}_1, \hat{\eta}_1, \hat{\eta}_2)$ and $\breve{\theta}=(\breve{\omega}, \breve{\lambda},\breve{\alpha}_1, \breve{\beta}_1, \breve{\eta}_1, \breve{\eta}_2)$ in $\Theta$, we have $g(\varepsilon_t;\hat{\theta}^{(1)})=g(\varepsilon_t;\breve{\theta}^{(1)}) \Rightarrow \hat{\omega}=\breve{\omega}$ and $c(\varepsilon_t;\hat{\theta}^{(2)}) =c(\varepsilon_t;\breve{\theta}^{(2)}) \Rightarrow \hat{\alpha}_1(|\varepsilon_{t-1} - \hat{\eta}_2 |-  \hat{\eta}_1 (\varepsilon_{t-1} - \hat{\eta}_2))^{\gamma} + \hat{\beta}_1 = \breve{\alpha}_1(|\varepsilon_{t-1} - \breve{\eta}_2 |-  \breve{\eta}_1 (\varepsilon_{t-1} - \breve{\eta}_2))^{\gamma} + \breve{\beta}_1$. 

We set $\varepsilon_t=\hat{\eta}_2$ and then $\varepsilon_t=\breve{\eta}_2$ and obtain the relation $\hat{\alpha}_1(|\breve{\eta}_2 - \hat{\eta}_2 |-  \hat{\eta}_1 (\breve{\eta}_2 - \hat{\eta}_2))^{\gamma} + \breve{\alpha}_1(|\hat{\eta}_2 - \breve{\eta}_2 |-  \breve{\eta}_1 (\hat{\eta}_2 - \breve{\eta}_2))^{\gamma}=0$. Since $\hat{\alpha}_1 >0$, $\breve{\alpha}_1>0$ and $|\breve{\eta}_1|<1, \ |\breve{\eta}_1|<1$ we should have $\hat{\eta}_2=\breve{\eta}_2$ and it follows $\hat{\beta}_1=\breve{\beta}_1$. Since $\varepsilon_t$ can take positive as well as negative values, we get $\hat{\eta}_1=\breve{\eta}_1$ and $\hat{\alpha}_1=\breve{\alpha}_1$.

$4.:$ There exist optimum points, namely $\theta_1^{*}=(\overline{\omega},\overline{\lambda},\overline{\alpha_1},\overline{\beta_1},\underline{\eta_1},\underline{\eta_2})$, $\theta_2^{*}=(\overline{\omega},\overline{\lambda},\overline{\alpha_1},\overline{\beta_1},\underline{\eta_1},\overline{\eta_2})$, $\theta_3^{*}=(\overline{\omega},\overline{\lambda},\overline{\alpha_1},\overline{\beta_1},\overline{\eta_1},\underline{\eta_2})$ and  $\theta_4^{*}=(\overline{\omega},\overline{\lambda},\overline{\alpha_1},\overline{\beta_1},\overline{\eta_1},\overline{\eta_2})$.

\subsubsection{Asymptotic Normality (Prop. \ref{prop: ass norm of QMLE})}

Since the fGARCH as defined above is not differentiable with respect to $\eta_2$, we cannot apply our theory developed for the asymptotic normality. So we will use a subfamily of this. By setting $\eta_2=0$ and $\gamma = \delta$ we obtain the family of \textit{asymmetric power ARCH/GARCH (apARCH)} models of \cite{DindGranEngl1993}, which subsumes many well known GARCH models (standard GARCH, absolute value GARCH, GJR GARCH, Threshold GARCH, non-linear ARCH,...). It reads
\begin{equation*}
R_t=\sigma_t \varepsilon_t, \ \ \ \sigma_t^{\delta}=\underbrace{\omega}_{g(\varepsilon_{t-1};\theta)}+ \underbrace{\lambda u_1(x_{t-1})}_{u(x_{t-1};\theta)}+ \underbrace{\left[\alpha_1(|\varepsilon_{t-1}|-  \eta_1 \varepsilon_{t-1})^{\delta} + \beta_1 \right]}_{c(\varepsilon_{t-1};\theta)} \sigma_{t-1}^{\delta} \ .
\label{eq:apARCH1}
\end{equation*} 
Since we are looking for the properties of the MLE given the observations $R_1,...,R_n$, these are considered as given. So we can rewrite the model and redefine the functions $c$ and $g$ as follows:
\begin{equation*}
R_t=\sigma_t \varepsilon_t, \ \ \ \sigma_t^{\delta}=\underbrace{\omega+ \alpha_1(|R_{t-1}|- \eta_1 R_{t-1})^{\delta}}_{g(\varepsilon_{t-1};\theta)} + \underbrace{\lambda u(x_{t-1})}_{u(x_{t-1};\theta)}+ \underbrace{ \beta_1 }_{c(\varepsilon_{t-1};\theta)} \sigma_{t-1}^{\delta} \ .
\label{eq:apARCH2}
\end{equation*} 

As we will see below this second representation has a big advantage. Let us now check that the assumptions \ref{ass: assumptions for the ass norm} are fulfilled, since these ensure the asymptotic normality of the likelihood estimators. 

$1.$ is clear.    

$2.$ is fulfilled due to the assumptions made on $\{\varepsilon_t\}_t$.

$3.$ holds since $g$ and $c$ and their first and second derivative exist and are continuous as we can see from the following equations
$$ \frac{\partial c(\varepsilon_t;\theta)}{\partial \theta} = \left(\frac{\partial c }{\partial\omega}, \frac{\partial c}{\partial\lambda}, \frac{\partial c}{\partial\alpha_1}, \frac{\partial c}{\partial\eta_1}, \frac{\partial c}{\partial\beta_1}\right)= (0, 0, 0, 0, 1) \ \ \textrm{ and } \ \  \frac{\partial^2 c(\varepsilon_t;\theta)}{\partial\theta\partial\theta^{'}} = 0_{5\times5} \ .$$
$$ \frac{\partial g(\varepsilon_t;\theta)}{\partial \theta} = \left(\frac{\partial g }{\partial\omega}, \frac{\partial g}{\partial\lambda}, \frac{\partial g}{\partial\alpha_1}, \frac{\partial g}{\partial\eta_1}, \frac{\partial g}{\partial\beta_1}\right)= \left(1, 0, (|R_t|-\eta_1 R_t)^{\delta}, -\alpha_1 \delta R_t(|R_t|-\eta_1 R_t)^{\delta-1}, 0 \right)$$
$$ \textrm{ and }\  \  \  \frac{\partial^2 g(\varepsilon_t;\theta)}{\partial\theta_i\partial\theta_j} = \left\{ 
\begin{array}{cl}
- \delta R_t(|R_t|-\eta_1 R_t)^{\delta-1} =: a_1  \ , & \textrm{if } (\theta_i,\theta_j) \in \{(\alpha_1,\eta_1), (\eta_1,\alpha_1)\} \\
\alpha_1 \delta (\delta-1) R_t^2(|R_t|-\eta_1 R_t)^{\delta-2} =: a_2 \ , &  \textrm{if } (\theta_i,\theta_j)=(\eta_1,\eta_1) \\
0 \ ,& \textrm{else.}
\end{array}
\right.$$

$4.:$ Since there exists (from proposition 3.1.) a $\alpha$ such that $R_t$ is $\alpha \delta$ stationary, then (using the Cauchy-Schwartz inequality) there exists a $\alpha^{'}\leq \alpha$ such that $E\left(|a_1|^{\alpha^{'}}+|a_2|^{\alpha^{'}}\right)<\infty.$ This implies $4.$ 

$5.:$ The first and second derivative of $c$ are trivially bounded on $\Theta$ since they are constant functions. Concerning the function $g$, all its constant derivatives are also trivially bounded. Let us look at the non constant derivatives: $|\partial g/\partial\alpha_1|$  and $|a_1|$ are bounded at $\underline{\eta_1}$ or $\overline{\eta_1}$, $|\partial g/\partial\eta_1|$ and $|a_2|$ are bounded at $(\overline{\alpha_1}, \underline{\eta_1})$ or $(\overline{\alpha_1}, \overline{\eta_1})$.

$6.$ and $7.$: Since 
\begin{eqnarray*}
\sigma_t^{\delta} &=& \sum\limits_{k=0}^{\infty} \left(\left[g(\varepsilon_{t-1-k}; \theta)+\lambda u_1(x_{t-1-k})\right] \prod\limits_{j=0}^{k-1} c(\varepsilon_{t-1-j}; \theta)\right) \\
&=& \sum\limits_{k=0}^{\infty} \left[\omega + \alpha_1(|R_{t-1}|- \eta_1 R_{t-1})^{\delta} +\lambda u_1(x_{t-1-k})\right] \beta^k  \ ,
\end{eqnarray*}
we get
\begin{eqnarray*}
\frac{1}{\sigma_t^{\delta}} \frac{\partial\sigma_t^{\delta}}{\partial\omega} &=& \frac{1}{\sigma_t^{\delta}} \sum\limits_{k=0}^{\infty} \beta^k \leq \frac{1}{\underline{\omega}} \sum\limits_{k=0}^{\infty} \overline{\beta}_1^k < \infty  \ , \\
\frac{\lambda}{\sigma_t^{\delta}} \frac{\partial \sigma_t^{\delta}}{\partial\lambda} &=& \frac{\lambda}{\sigma_t^{\delta}} \sum\limits_{k=0}^{\infty} u(x_{t-1-k}) \beta_1^k  = \frac{1}{\sigma_t^{\delta}} \sum\limits_{k=0}^{\infty} \lambda u(x_{t-1-k}) \beta_1^k \leq \frac{\sigma_t^{\delta}}{\sigma_t^{\delta}} =1 \  \  \Longrightarrow \  \  \frac{1}{\sigma_t^{\delta}} \frac{\partial \sigma_t^{\delta}} {\partial\lambda} \leq \frac{1}{\underline{\lambda}} \ , \\
\frac{\alpha_1}{\sigma_t^{\delta}} \frac{\partial \sigma_t^{\delta}}{\partial\alpha_1} &=& \frac{\alpha_1}{\sigma_t^{\delta}} \sum\limits_{k=0}^{\infty} (|R_{t-1-k}|- \eta_1 R_{t-1-k})^{\delta} \beta_1^k \leq \frac{\sigma_t^{\delta}}{\sigma_t^{\delta}} =1 \  \  \Longrightarrow \  \  \frac{1}{\sigma_t^{\delta}} \frac{\partial \sigma_t^{\delta}} {\partial\alpha_1} \leq \frac{1}{\underline{\alpha_1}} \ ,  \\
J_{0,0}^{-1} \frac{\partial \sigma_t^{\delta}}{\partial\eta_1} &=& J_{0,0}^{-1} \frac{1}{\sigma_t^{\delta}} \left| \sum\limits_{k=0}^{\infty} -\alpha_1 R_{t-1} (|R_{t-1}|- \eta_1 R_{t-1})^{\delta-1} \beta_1^k \right|   
\leq \frac{\sigma_t^{\delta}}{\sigma_t^{\delta}} \ \ \Longrightarrow \ \ \frac{1}{\sigma_t^{\delta}} \frac{\partial \sigma_t^{\delta}} {\partial\eta_1} \leq J_{0,0} \ ,
\end{eqnarray*}
where $J_{0,0}:= \left(\min\{|1+\underline{\eta}_1|,|-1+\overline{\eta}_1|\} / \delta \right)^{-1}$.

Finally, due to the fact for all $x\geq0$ and all $0< s <\alpha$, $x/(1+x) \leq x^{s}$, we have
\begin{eqnarray}
E_{\theta_0} \left(\sup\limits_{\theta\in\Theta} \frac{\beta_1}{\sigma_t^{\delta}} \frac{\partial \sigma_t^{\delta}}{\partial\beta}\right) &=&  E_{\theta_0} \sup\limits_{\theta\in\Theta}  \frac{\beta_1}{\sigma_t^{\delta}} \sum\limits_{k=1}^{\infty} [g(\varepsilon_{t-1-k})+ u(x_{t-1-k})] k \beta_1^{k-1} \nonumber \\
&\leq& E_{\theta_0} \sup\limits_{\theta\in\Theta}  \sum\limits_{k=1}^{\infty} k \frac{[g(\varepsilon_{t-1-k})+ u(x_{t-1-k})] \beta_1^k} {\omega+[g(\varepsilon_{t-1-k})+ u(x_{t-1-k})] \beta_1^k} \nonumber \\
&\leq& E_{\theta_0} \sup\limits_{\theta\in\Theta}  \sum\limits_{k=1}^{\infty} k \left(\frac{[g(\varepsilon_{t-1-k})+ u(x_{t-1-k})]} {\omega}\right)^{s}\beta_1^{s k} \nonumber \\
&\leq& \frac{1}{{\underline{\omega}}} \sum\limits_{k=1}^{\infty} k \overline{\beta}_1^{s k} \underbrace{ E_{\theta_0} \sup\limits_{\theta\in\Theta} \left[g(\varepsilon_{t-1-k})+ u(x_{t-1-k}) \right]^{s}}_{=:J_{1} <\infty \textrm{ since optimum points for $g$ and $u$ exist}} \nonumber \\
&=& \frac{J_{1}}{\underline{\omega}} \sum\limits_{k=1}^{\infty} k \overline{\beta}_1^{s k} =: J_{1,1} < \infty.
\label{eq: eq1}
\end{eqnarray}
Proceeding as above (see equation \eqref{eq: eq1}) and using the relation $x/(1+x) \leq x^{s/2}$ we get 
\begin{eqnarray*}
E_{\theta_0} \left(\sup\limits_{\theta\in\Theta} \frac{1}{\sigma_t^{\delta}} \frac{\partial \sigma_t^{\delta}}{\partial\beta}\right)^2 \leq \frac{J_{2,2}}{\underline{\beta}_1^2} < \infty \ \ \textrm{with} \ \ J_{2,2} < \infty.
\end{eqnarray*}
Using the Cauchy-Schwartz inequality over the derivative with respect to $\beta_1$ and the boundedness of the other derivative we get that $A=E_{\theta_0}\left(\frac{1}{\sigma_t^{2\delta}(\theta_0)} \frac{\partial\sigma_t^{\delta}(\theta_0)}{\partial \theta} \frac{\partial\sigma_t^{\delta}(\theta_0)}{\partial \theta^{'}}\right)$ exists. We first remark that $A$ is positive semi definite because of its form. Let us now suppose that there exists a vector $m=(m_1,...,m_5)$ such that $0=m^{'} \cdot \frac{\partial\sigma_{t+1}^{\delta}(\theta_0)}{\partial \theta}$. This should -due to stationarity of the derivative- hold at all times  $t$. Since $ \ \frac{\partial\sigma_{t+1}^{\delta}(\theta_0)}{\partial \theta}= (1, a_3:=(|R_t|- \eta_1 R_t)^{\delta}, a_4:= -\alpha_1 \delta R_t(|R_t|- \eta_1 R_t)^{\delta-1}, u(x_t), \sigma_t^{\delta} )^{'} + \frac{\partial\sigma_t^{\delta}(\theta_0)}{\partial \theta} \ $ we should have $m_1 + m_2 a_3 + m_3 a_4 + m_4 u(x_t) + m_5 \sigma_t^{\delta} =0$. Due to the fact that $x_t$ is independent of $\varepsilon_t$ and $x_{t-1},x_{t-2},...$ we should have $m_4=0$. Furthermore $a_3$ and $a_4$ are multiples of $\varepsilon_t$ which is non-degenerate and independent of $\sigma_t^{\delta}$, thus  $m_2 a_3 + m_3 a_4=0$, i.e. $(|R_t|- \eta_1 R_t)^{\delta-1} \left[m_2(|R_t|- \eta_1 R_t)-\alpha_1 \delta m_3 R_t\right]=0\ $ Since $\varepsilon_t$ is non-degenerate and $\mathbb{P}(\sigma_t>0)=1$ we should have $(m_2(|R_t|- \eta_1 R_t)-\alpha_1 \delta m_3 R_t)=0$. Using the same arguments together with the fact that $\alpha_1 \delta >0$ we have $m_2=m_3=0$. $m_5=0$ is due to the fact that $\sigma_t$ is non constant. It follows $m_1=0$ and thus that $A$ is indeed positive definite. 

Concerning the second derivative we have:
\begin{eqnarray*}
\frac{\partial^2 \sigma_t^{\delta}}{\partial\omega^2}= \frac{\partial^2 \sigma_t^{\delta}}{\partial\omega \partial\lambda}= \frac{\partial^2 \sigma_t^{\delta}}{\partial\omega \partial\alpha_1}= \frac{\partial^2 \sigma_t^{\delta}}{\partial\omega \partial\eta_1}=  \frac{\partial^2 \sigma_t^{\delta}}{\partial\lambda^2}= \frac{\partial^2 \sigma_t^{\delta}}{\partial\lambda \partial\alpha_1}= \frac{\partial^2 \sigma_t^{\delta}}{\partial\lambda \partial\eta_1}= \frac{\partial^2 \sigma_t^{\delta}}{\partial\alpha_1^2}=0 \ ,
\end{eqnarray*}
\begin{eqnarray*}
\frac{1}{\sigma_t^{\delta}}\frac{\partial^2 \sigma_t^{\delta}}{\partial\beta_1 \partial\omega} = \frac{1}{\sigma_t^{\delta}}\sum\limits_{k=1}^{\infty} k \beta_1^{k-1} \leq \frac{1}{\underline{\omega}} \sum\limits_{k=1}^{\infty} k \overline{\beta}_1^{k-1} <\infty, \ \frac{1}{\sigma_t^{\delta}}\frac{\partial^2 \sigma_t^{\delta}}{\partial\eta_1 \partial\alpha_1} \leq \frac{1}{\underline{\alpha}_1}J_{0,0} \ , \  \  \frac{1}{\sigma_t^{\delta}}\frac{\partial^2 \sigma_t^{\delta}}{\partial\eta_1^2} \leq \frac{\delta-1}{\delta} J_{0,0}^2 \ .
\end{eqnarray*}
Moreover using the same technique as in equation \eqref{eq: eq1} we have for all $s<\alpha$:
\begin{eqnarray*}
E_{\theta_0} \left(\sup\limits_{\theta\in\Theta}\frac{\partial^2 \sigma_t^{\delta}}{\partial\beta_1 \partial\lambda}\right) \leq \frac{J_{3,3}}{\underline{\beta}_1}, \ \ E_{\theta_0} \left(\sup\limits_{\theta\in\Theta}\frac{\partial^2 \sigma_t^{\delta}}{\partial\beta_1 \partial\eta_1}\right) \leq \frac{J_{4,4}}{\underline{\beta}_1}, \ \ E_{\theta_0} \left(\sup\limits_{\theta\in\Theta}\frac{\partial^2 \sigma_t^{\delta}}{\partial\beta_1 \partial\alpha_1}\right) \leq \frac{J_{4,4}}{\underline{\beta}_1} J_{0,0}\\
E_{\theta_0} \left(\sup\limits_{\theta\in\Theta}\frac{\partial^2 \sigma_t^{\delta}}{\partial\beta_1 \partial\lambda}\right) \leq \frac{J_{5,5}}{\underline{\beta}_1^2}, \ \
\end{eqnarray*}
where $J_{3,3}, J_{4,4}, J_{5,5}, J_{6,6} <\infty$.

The same technique is also used for the third derivative of $\sigma_t^{\delta}$. In this way we get $6.$
\section{Applications}
In this section we give some applications of the properties that we presented and proved before, in order to show some useful consequences of these properties. These applications are by far not all applications that can be derived.

\subsection{Use of ergodicity for value at risk calculation}
The value at risk (VaR) is nowadays one of the most important concepts in risk management. It is prescribed by regulatory requirements (e.g. in Basel II) and also used internally by companies for risk management purposes. The value at risk over a portfolio answers the question: what is the maximum loss on the portfolio with a given probability (level) over a given time horizon? In Basel II for example the level is 99\% and the horizon is fixed to ten working days (i.e. two weeks). Internally the companies are mostly interested in the VaR at level $\alpha = 95$\%.

Mathematically the VaR can be defined as:  
\begin{eqnarray*}
VaR^{\alpha}_{t,t+h}=\inf\left\{x\vert P(r_{t+h}\leq x\vert\mathcal{F}_{t})\geq 1-\alpha\right\},
\end{eqnarray*}
where $\mathcal{F}_{t}$ denotes the set of all information which is available until time $t$, $r_{t+h}=\ln(P_{t+h}/P_t)$ denotes the $h-$periodic log-return at time $t+h$ and $P_t$ represents the portfolio value at time $t$.

In this formula we get a value at risk in terms of log-returns. We can easily convert it to a value at risk in terms of returns $\left(\exp(VaR^{\alpha}_{t,t+h})-1\right)$ or in terms of portfolio value $\left(P_t \cdot \left[\exp(VaR^{\alpha}_{t,t+h})-1\right]\right)$.

Since our ergodicity result is for univariate GARCH models, we can use it in the following way: If the portfolio contains just one asset we directly model our asset through the GARCH. If in contrary the portfolio contains many assets, we directly model the portfolio (as one security) and not the individual securities\footnote{In order to model the securities individually we would need a multidimensional GARCH. However our ergodicity result is designed for univariate GARCH models.}. 

We suppose that we are at time $t$ and want to compute the VaR in $h$ periods i.e. $VaR^{\alpha}_{t,t+h}$ (e.g. for one year, $h= 250$ working days, i.e. we compute $VaR^{\alpha}_{t,t+250}$). The first method is a usual Monte Carlo algorithm. In the usual Monte Carlo algorithm we have to simulate many independent paths, says $n$ paths (for example $n=40000$). Since the simulation of each path requires to draw $2 \times h \ (= 2 \times 250 = 500)$ random numbers, namely the process $\{(\varepsilon_t, x_t)\}_t$, we need in total $2 \times h \times n (= 2 \times 250 \times 40000 = 2 \cdot 10^7)$ random numbers to compute the  $VaR^{\alpha}_{t,t+h}$.

If we now use the ergodic property, we can, instead of simulating many independent paths, work with just one long path of length $h-1+n \ (=249+40000)$, since the ergodicity property implies that the VaR computed through independent paths is nearly equal to the one computed with one path, given that the process is ergodic. As we mentioned in the third section the real process which started at some time $0$ with a predefined volatility's distribution is not ergodic, however converges with an exponential rate towards an ergodic process. Thus by using just one long path we need a burn-in period; this means we do not consider or take into account the first $N_b$ returns $r_t$ and volatility $\sigma_t$ that we will simulate. After the time $N_b$ the processes $r_t$ and $\sigma_t$ are considered as stationary. We can for example take $N_b=1000$ as length of the burn-in period. Thus we need $2\times(N_b+h-1+n) \ (=2(1000+249+40000)=82498)$ random numbers\footnote{Note that we have multiplied by $2$, since at each time step we need to simulate two random numbers, $\{\varepsilon_t\}_t$ and $\{x_t\}_t$.}.

Comparing the two methods we get the following ratio:
\begin{equation}
\frac{\mbox{\scriptsize Number of simulations for method 1}}{\mbox{\scriptsize Number of simulations for method 2}}=\frac{2\cdot h \cdot n}{2(N_b+h-1+n)} \approx h \ (\mbox{for large n})\ ,
\end{equation}
where ``method 1'' denotes the usual Monte Carlo algorithm with independently simulated paths and ``method 2'' is the one in which we use the ergodic property. 

Hence the ergodic property can help to reduce the computational time and eventually the required memory.

We now want to illustrate that both VaRs computed through the two different methods are almost identically by the following example: we consider the $GJR-GARCH(1,1)$ of \cite{GlostJagaRunk1993}. We assume the noise process $\{\varepsilon_t\}_t$ is $iid \ \mathcal{N}(0,1)$. For the exogenous variable we assume an autoregressive process of order $1$ ($AR(1)$), i.e. $x_t=\phi x_{t-1} + \eta_t$ where the noise process $\{\eta_t\}_t$ is $iid \ Cauchy(1)$. The autoregressive component $\phi$ equals $0.8$, which implies that $\{x_t\}_t$ is stationary and ergodic. For the transformation function $u$ of the exogenous process we consider the absolute value function, $u(x_t)=|x_t|$. The processes  $\{\eta_t\}_t$ and $\{\varepsilon_t\}_t$ are independent so that the assumption \ref{assumption on noise processes} is fulfilled. Our model thus reads
\begin{eqnarray*}
R_t&=&\sigma_t \varepsilon_t \ , \\
\sigma_t^{2}&=& \omega + \alpha_1 R_{t-1}^2 + \beta_1 \sigma_{t-1}^2 + \gamma_1 R_{t-1}^2 1_{\{R_{t-1}<0\}} + \lambda |x_{t-1}|  \ . 
\end{eqnarray*}
With the following parameters\footnote{The  parameters $h$ and $\alpha$ are chosen in the following way: $h=10$ (10 workings days or two weeks) and $\alpha=99\%$, so that they meet the regulatory requirements. The parameter $N_{ex}$ represents the burn-in period needed in order to obtain a stationary sample from the exogenous AR(1) process $\{x_t\}_t$.} $\omega=0.04, \ \alpha_1=0.1, \ \beta_1=0.8, \ \gamma_1=0.06, \ \lambda=0.02, \ \sigma_0=0.02, \ r_0=0, \ h=10, \ n=100000, \ N_b=5000, \ N_{ex}=10000$, we compute the values at risk several times ($500$ times) at level $\alpha=99$\% with both methods (namely one long path and many independent paths). Conducting an approximative $t-test$\footnote{We conduce a two-sample two-sided t-test to compare the means of the results obtain from both methods.} to compare the results of both methods we obtain as $p-$value $p=0.3609$, which indicates that both results are very similar.

\subsection{Computation of the confidence region}

The asymptotic normality property proven in the fifth section can be used to compute confidence regions of the true value. This region is an \textit{asymptotic} region since we proved \textit{asymptotic} normality of the likelihood estimators.

Since we proved that $\sqrt{n}(\hat{\theta}_n-\theta_0) \ \longrightarrow  \ \mathcal{N}\left(0,B\right)$, with $B= \frac{\delta^2}{4}(\kappa-1) \left( E_{\theta_0} \left(\frac{1}{\sigma_t^{2\delta}(\theta_0)} \frac{\partial\sigma_t^{\delta}(\theta_0)}{\partial \theta} \frac{\partial\sigma_t^{\delta} (\theta_0)}{\partial \theta^{'}}\right) \right)^{-1}$, we get that the ($1-p$)-confidence region\footnote{Typically $p$ equals 5\% or 1\%.} for $\theta$ is the set
\begin{equation}
CR_n(1-p) \approx \left\{ \theta : \  n(\theta-\hat{\theta}_n)^{'}B_n^{-1}(\theta-\hat{\theta}_n) \leq \chi_m^2(1-p) \right\}  \ ,
\label{eq: confidence bound}
\end{equation}
where $\chi_m^2(1-p)$ denotes the $1-p$ quantile of a chi-squared distribution with $m$ degrees of freedom, $m$ denotes the number of elements of $\theta$, $B_n:=\frac{\delta^2}{4}(\kappa-1) \left( E_{\theta_0} \left(\frac{1}{\sigma_t^{2\delta}(\hat{\theta}_n)} \frac{\partial\sigma_t^{\delta}(\hat{\theta}_n)}{\partial \theta} \frac{\partial\sigma_t^{\delta}(\hat{\theta}_n)} {\partial \theta^{'}}\right) \right)^{-1}$ converges in probability towards $B$.

If we are only interested in some parameters $\theta_i \ i \in I$ of $\theta$, the same result \eqref{eq: confidence bound} applies, however with a sub-matrix $B=(B_{i,j})_{i,j\in I}$ containing only the parameters under interest and $m=|I|$\footnote{$|I|$ denotes the cardinal of $I$, i.e., the number of elements of $I$.}. We also note that if $|I|=1$ the confidence region is equivalent to the classical confidence interval that we usually obtain with the quantile of the standard normal distribution.      
\section{Conclusion and Outlook}
In this work we have studied a new type of GARCH models which is gaining importance in both theory and practice: the so-called GARCHX models. The GARCHX models are GARCH models augmented with an exogenous variable. We have studied a very general model which includes almost all well-known and used GARCH models. As exogenous variable we have assumed a function of a stationary and exogenous process. For these GARCHX models we found conditions for some characteristics of the model, namely the stationarity, ergodicity, geometric ergodicity and existence of moments. We have pointed out and solved some practical aspects related to some properties listed above. After this we also checked one very important aspect of the model when it is fitted to data, namely the asymptotic theory. We have given the conditions for the quasi maximum likelihood estimator to be consistent and asymptotically normal distributed. After this, an example concerning a whole family of GARCH models augmented with an exogenous variable was studied: we have shown that under mild conditions this family (or a sub-family) fulfills the requirements needed to get the above enunciated properties. Finally we have given some practical applications of the properties that we studied: we demonstrated how we can use the stationarity and ergodicity properties in order to quickly compute the value at risk. We also studied the use of the asymptotic normality property in order to find the confidence region. 

To extend this work we can look at the same model family, however with the exogenous variable being non stationary as in this paper. A good starting point of such a research would certainly be the work done by \cite{HanKris2012} who studied asymptotic properties within the standard GARCH model and as exogenous variable, a quadratic function of a fractionally integrated process. 
\section{Appendix}

\textbf{\textit{Proof of the proposition \ref{propoStatErgoInfty}}}

We set for all $n \in \N $ $$S_n:=g(\varepsilon_{t-1})+u(x_{t-1}) + \sum\limits_{k=1}^{n}\prod_{j=0}^{k-1} c(\varepsilon_{t-1-j}) \left[g(\varepsilon_{t-1-k}) + u(x_{t-1-k})\right] \ .$$

In order to write short we consider that $\prod\limits_{j=0}^{-1} c(\varepsilon_{t-1-j})=1$, then we can shortly write:
\begin{equation}
S_n:=\sum\limits_{k=0}^{n}\prod_{j=0}^{k-1} c(\varepsilon_{t-1-j}) \left[g(\varepsilon_{t-1-k}) + u(x_{t-1-k})\right] \ .
\end{equation} 
It follows
\begin{eqnarray*}
E\left(\lim\limits_{n\rightarrow \infty} S_n \right)^{\alpha} &=&E\left(\lim\limits_{n\rightarrow \infty} \sum\limits_{k=0}^{n}\prod_{j=0}^{k-1} c(\varepsilon_{t-1-j}) \left[g(\varepsilon_{t-1-k}) + u(x_{t-1-k})\right] \right)^{\alpha} \\
&\leq& E\underbrace{\sum\limits_{k=0}^{\infty}\prod_{j=0}^{k-1} \left[c(\varepsilon_{t-1-j}) g(\varepsilon_{t-1-k}) \right]^{\alpha}}_{A} 
 + E\underbrace{\sum\limits_{k=0}^{\infty}\prod_{j=0}^{k-1} \left[c(\varepsilon_{t-1-j}) u(x_{t-1-k}) \right]^{\alpha}}_{B}  \ ,
\end{eqnarray*}
where the inequality is due to the fact that $\alpha \in (0,1] $ from which follows that $(a+b)^{\alpha}\leq a^{\alpha}+b^{\alpha}$ for all $a,b>0$.
\begin{eqnarray*}
E(A) = E \left[\sum\limits_{k=0}^{\infty}\prod_{j=0}^{k-1} \left[c(\varepsilon_{t-1-j})\right]^{\alpha} \left[g(\varepsilon_{t-1-k}) \right]^{\alpha} \right] 
 &=& \sum\limits_{k=0}^{\infty} E \left[g(\varepsilon_{t-1-k}) \right]^{\alpha} E\left[\prod_{j=0}^{k-1} c(\varepsilon_{t-1-j})\right]^{\alpha} \\
 &=& E \left[g(\varepsilon_{t}) \right]^{\alpha} \sum\limits_{k=0}^{\infty}  \left[ E\left( c(\varepsilon_{t})\right)^{\alpha}\right]^{k} \\
 &=& E \left[g(\varepsilon_{t}) \right]^{\alpha} \left[1-E\left[c(\varepsilon_t)\right]^{\alpha}\right]^{-1}\\
 &<& \infty \ .
\end{eqnarray*}
Note that the second and the third equalities are implied from the assumption \ref{assumption on noise processes} which states that the process $\{(\varepsilon_t, x_t)\}_t$ is $i.i.d.$. The last equality is due to the fact that $E\left[c(\varepsilon_t)\right]^{\alpha} < 1$.

Similarly,
\begin{eqnarray*}
E(B) = E \left[\sum\limits_{k=0}^{\infty}\prod_{j=0}^{k-1} \left[c(\varepsilon_{t-1-j})\right]^{\alpha} \left[u(x_{t-1-k}) \right]^{\alpha} \right] 
 &=& \sum\limits_{k=0}^{\infty} E \left[u(x_{t-1-k}) \right]^{\alpha} E\left[\prod_{j=0}^{k-1} c(\varepsilon_{t-1-j})\right]^{\alpha} \\
 &=& E \left[u(x_{t}) \right]^{\alpha} \sum\limits_{k=0}^{\infty}  E\left[\left( c(\varepsilon_{t-1-j})\right)^{\alpha}\right]^{k} \\
 &=& E \left[u(x_{t}) \right]^{\alpha} \left[1-E\left[c(\varepsilon_t)\right]^{\alpha}\right]^{-1}\\
 &<& \infty \ .
\end{eqnarray*}
The second line follows from the fact that $x_{t-1-k}$ is independent from $\varepsilon_{t-1-j}$, $j=0,...,k-1$.

$E(A)+E(B)<\infty$ implies that the sequence $(S_n)$ converges almost surely (a.s.) through a limit that we denote by $\sigma_t^{\delta}$, this means 
\begin{equation*}
\sigma_t^{\delta}=\sum\limits_{k=0}^{\infty} \prod\limits_{j=0}^{k-1} c(\varepsilon_{t-1-j})\left[g(\varepsilon_{t-1-k}) + u(x_{t-1-k})\right]  \ .
\end{equation*}
Thus $R_t=\sigma_t \varepsilon_t$ is a fixed functional of iid random variables $(\varepsilon_j,x_j)$, so it is stationary and ergodic. 

Furthermore $E\left|R_t\right|^{\alpha \delta}=E\left|\sigma_t \varepsilon_t \right|^{\alpha \delta} = \underbrace{E\left|\sigma_t \right|^{\alpha \delta}}_{<\infty} \cdot \underbrace{ E\left| \varepsilon_t \right|^{\alpha \delta} }_{<\infty}$. It follows that the process $\{R_t\}_t$ is  $\alpha \delta$-stationary.

The uniqueness is quickly proved by assuming that there is another solution for our model, i.e. $R_t=h_t \varepsilon_t$ with $h_t$ solving the equation \eqref{eq:second equat of the model}, then
\begin{eqnarray*}
E\left|h_t^{\delta} - \sigma_t^{\delta} \right|^{\alpha} = E\left(c(\varepsilon_{t-1})\right)^{\alpha} E\left|h_{t-1}^{\delta} - \sigma_{t-1}^{\delta} \right|^{\alpha} 
 &=&  \left[E\left(c(\varepsilon_{t-1})\right)^{\alpha}\right]^2 E\left|h_{t-2}^{\delta} - \sigma_{t-2}^{\delta} \right|^{\alpha} \\
& \vdots & \\
& = & \left[E\left(c(\varepsilon_{t-1})\right)^{\alpha}\right]^k E\left|h_{t-k}^{\delta} - \sigma_{t-k}^{\delta} \right|^{\alpha} \\
& \leq & \left[\underbrace{E\left(c(\varepsilon_{t-1})\right)^{\alpha}}_{<1}\right]^k \underbrace{\left( E\left|h_{t-k}^{\delta}\right|^{\alpha} + E\left|\sigma_{t-k}^{\delta} \right|^{\alpha} \right)}_{<\infty} \\
&& \longrightarrow 0 \textrm{ as } k \longrightarrow \infty  \ .    \hspace{3cm} \square
\end{eqnarray*}

\textbf{\textit{Proof of the proposition \ref{propoStatErgoTime0}}}

To prove equation \eqref{eq1:propoStatErgoTime0} we first prove by recursion (over $t \in \N $) the following equality:
$ \tilde{\sigma}_t^{\delta} =\sigma_t^{\delta} + \xi_t$ where $\xi_t= (\tilde{\sigma}_0^{\delta} -\sigma_0^{\delta}) \prod\limits_{j=0}^{t-1}c(\varepsilon_j)$ and we set $\prod\limits_{j=0}^{-1}c(\varepsilon_j)$ equals 1. 
\begin{itemize}
	\item for $t=0$ we have $\xi_0=\tilde{\sigma}_0^{\delta} -\sigma_0^{\delta}$  by assumption.
	\item $t-1 \longrightarrow t$: we suppose that the assertion is fulfilled at time $t-1$.
	   \begin{eqnarray*}
     \tilde{\sigma}_t^{\delta}&=& g(\varepsilon_{t-1}) + u(x_{t-1}) + c(\varepsilon_{t-1})\tilde{\sigma}_{t-1}^{\delta} \\
                      &=& g(\varepsilon_{t-1}) + u(x_{t-1}) + c(\varepsilon_{t-1})\left[\sigma_{t-1}^{\delta} + \xi_{t-1} \right] \\
                      &=& g(\varepsilon_{t-1})+u(x_{t-1})+c(\varepsilon_{t-1})\sigma_{t-1}^{\delta}+ c(\varepsilon_{t-1})\cdot ( \tilde{\sigma}_0^{\delta} -\sigma_0^{\delta}) \prod\limits_{j=0}^{t-2}c(\varepsilon_j)\\
                      &=& \sigma_t^{\delta} + \xi_t  \ .
     \end{eqnarray*} 
\end{itemize}
It follows
\begin{eqnarray*}
E\left|\xi_t\right|^{\alpha} = E\left[\left|\tilde{\sigma}_0^{\delta}-\sigma_0^{\delta}\right|^{\alpha} \prod\limits_{j=0}^{t-1}(c(\varepsilon_j))^{\alpha}\right] 
  &=& E\left|\tilde{\sigma}_0^{\delta}-\sigma_0^{\delta}\right|^{\alpha}\prod\limits_{j=0}^{t-1}E(c(\varepsilon_j))^{\alpha} \\
  &\leq& E(\left|\tilde{\sigma}_0\right|^{\delta\alpha}+\left|\sigma_0\right|^{\delta \alpha}) \left(E[c(\varepsilon_1)]^{\alpha}\right)^t \\
  &=& (\underbrace{E\left|\tilde{\sigma}_0\right|^{\delta\alpha}}_{A_1}+ \underbrace{E\left|\sigma_0\right|^{\delta \alpha}}_{A_2}) \left(\underbrace{E[c(\varepsilon_1)]^{\alpha}}_{<1}\right)^t \\
  &=&\mathcal{O}(\rho^t) \ \ \textrm{where} \ \ \rho=E[c(\varepsilon_1)]^{\alpha},
\end{eqnarray*}
since $A_1 < \infty$ by assumption, $A_2 < \infty$ by proposition \ref{propoStatErgoInfty}. The above first inequality is due to the triangular inequality and the fact that $\alpha < 1$. This proves the equation \eqref{eq1:propoStatErgoTime0}.

From this we get the proof of the assertion \eqref{eq2:propoStatErgoTime0}:
\begin{eqnarray*}
\tilde{R}_t^{\delta}&=& \tilde{\sigma}_t^{\delta} \varepsilon_t^{\delta} 
    = (\sigma_t^{\delta} + \xi_t)\cdot \varepsilon_t^{\delta} 
    = \underbrace{\sigma_t^{\delta} \cdot \varepsilon_t^{\delta}}_{=R_t^{\delta}} + \underbrace{\xi_t \cdot \varepsilon_t^{\delta}}_{=:\psi_t} \\
\end{eqnarray*}
and 
\begin{equation*}
E\left|\psi_t\right|^{\alpha} = E\left|\varepsilon_t\right|^{\delta\alpha} E\left|\xi_t\right|^{\alpha} = \underbrace{E\left|\varepsilon_t\right|^{\delta\alpha}}_{<\infty} \cdot \mathcal{O}(\rho^t)= \mathcal{O}(\rho^t),
\end{equation*}
where the first equality of the above equation is due to the independence of $\xi_t$ and $\varepsilon_t$. $\hspace{1cm} \square$
\newline

\textbf{\textit{Proof of the proposition \ref{existence of moments}}}

$(i)\Rightarrow(ii):$ 
\begin{eqnarray*}
\infty > E\left|R_t\right|^{m\delta} &=& E \left|\varepsilon_t\right|^{m\delta} E \left|\sigma_t^{\delta}\right|^{m} \\
 &=& E \left|\varepsilon_t\right|^{m\delta}E\left[ \sum\limits_{k=0}^{\infty} \left[u(x_{t-1-k}) + g(\varepsilon_{t-1-k}) \right] \prod_{j=0}^{k-1} c(\varepsilon_{t-1-j}) \right]^m \\
 & \underbrace{\geq}_{m \geq 1}& E \left|\varepsilon_t\right|^{m\delta}E \sum\limits_{k=0}^{\infty} \left\{\left[u(x_{t-1-k})\right]^m + \left[g(\varepsilon_{t-1-k}) \right]^m \right\} \left[\prod_{j=0}^{k-1} c(\varepsilon_{t-1-j})\right]^m \\
 & =& E \left|\varepsilon_t\right|^{m\delta} \sum\limits_{k=0}^{\infty} \left\{E[g(\varepsilon_{t-1-k})]^m +E[u(x_{t-1-k})]^m\right\} \prod\limits_{j=0}^{k-1} E[c(\varepsilon_{t-1-j})]^m  \\
 &=& \underbrace{E\left|\varepsilon_t\right|^{m\delta}}_{<\infty} \sum\limits_{k=0}^{\infty} \left\{E[g(\varepsilon_{t})]^m +E[u(x_{t})]^m\right\} \left[E[c(\varepsilon_{t})]^m\right]^k\\
 &=& \underbrace{E\left|\varepsilon_t\right|^{m\delta}}_{<\infty} \underbrace{\left\{E[g(\varepsilon_{t})]^m +E[u(x_{t})]^m\right\}}_{<\infty} \sum\limits_{k=0}^{\infty} \left[E[c(\varepsilon_{t})]^m\right]^k.
\end{eqnarray*}
It implies $E[c(\varepsilon_{t})]^m < 1$.

Note that the inequality (``$\geq$'') is obtained from the binomial formula.

$(ii)\Rightarrow(i):$ 
we prove it by using the following criterion of \cite{Twee1988} over the variance process $\{\sigma_t\}_t$:
\begin{lemma}
Let $\{X_n\}_n$ be a weak Feller process over $(S,\mathcal{B})$, where $S$ is a locally compact separable metric space with $\mathcal{B}$ the Borel $\sigma$-field on $S$ and P(.,.) denoting the transition probability. We suppose that there exists a compact set $A \in \mathcal{B}$ such that
\begin{itemize}
	\item there exists  a non negative function $f$ and a real number $\varepsilon$ satisfying
	    \begin{equation}
	      \int_{A^c}P(x,dy)f(y) \leq f(x) - \varepsilon \ \ x \in A^c  \ ,
      \label{eq:lemma Tweedie 1}
      \end{equation} 
      then there exists a $\sigma$-finite invariant measure $\mu$ for $P$ with $0 < \mu(A) < \infty$.
  \item If
      \begin{equation}
      \int_{A}\mu(dx)\left[\int_{A^c}P(x,dy)f(y)\right] <\infty  \ ,
      \label{eq:lemma Tweedie 2}
      \end{equation}
      then $\mu$ is finite, and hence $\pi=\mu / \mu(S)$ is an invariant probability measure.
  \item If
      \begin{equation}
      \int_{A^c}P(x,dy)f(y) \leq f(x) - g(x), \ \ x \in A^c  \ ,
      \label{eq:lemma Tweedie 3}
      \end{equation}
      then $\mu$ admits a finite g-moment, i.e. $\int_S\mu(dy)g(y)<\infty.$
\end{itemize}
\label{lemma of Tweedie}
\end{lemma}
$\{\sigma_t^{\delta}\}_t$ is a Feller process, since it is a Markov process and for every bounded continuous function $f$ on $\mathbb{R}_{>0}$, the function $E\left[f(\sigma_t^{\delta})|\sigma_{t-1}^{\delta}=x\right]$ is continuous in $x$.

Furthermore for the function $f$ defined on $\mathbb{R}_{>0}$ by $f(x)=1+x^m$, we have for all $x>0$:
\begin{eqnarray*}
E\left[f(\sigma_t^{\delta})|\sigma_{t-1}^{\delta}=x\right] &=& 1 + E[g(\varepsilon_{t-1})+u(x_{t-1})+c(\varepsilon_{t-1})x]^m \\
&=& 1 + \sum\limits_{\stackrel{0 \leq i,j,k \leq m}{i+j+k=m}}E\left\{\left(\stackrel{m}{i,j,k}\right) (g(\varepsilon_{t-1}))^i (u(x_{t-1}))^j (c(\varepsilon_{t-1})x)^k \right\} \\
&=& 1 + E[c(\varepsilon_{t-1})]^m x^m \\
&& + \sum\limits_{\stackrel{0 \leq i,j,k \leq m}{i+j+k=m \ k \neq m}}E\left\{\left(\stackrel{m}{i,j,k}\right) (g(\varepsilon_{t-1}))^i (u(x_{t-1}))^j (c(\varepsilon_{t-1})x)^k \right\} \\
& \leq &  1 + E[c(\varepsilon_{t-1})]^m x^m + (1+x)^{m-1} \\
&& \times \underbrace{\sum\limits_{\stackrel{0 \leq i,j,k \leq m}{i+j+k=m \ k \neq m}}\left(\stackrel{m}{i,j,k}\right) \underbrace{ E\left\{(g(\varepsilon_{t-1}))^i (u(x_{t-1}))^j (c(\varepsilon_{t-1}))^k \right\}}_{=:A_{i,j,k}}}_{=:B}  \ .
\end{eqnarray*}
Note that the above inequality is due to the fact that $x^k \leq (1+x)^{m-1}$ for all $k\in{0,...,m-1}$.

Since $(g(\varepsilon_{t-1}))^i \in L^{\frac{m}{i}}(0,\infty), \ \ (u(x_{t-1}))^j \in L^{\frac{m}{j}}(0,\infty), \ \ (c(\varepsilon_{t-1}))^k \in L^{\frac{m}{k}}(0,\infty), \ \ \frac{1}{\frac{m}{i}}+ \frac{1}{\frac{m}{j}} +\frac{1}{\frac{m}{k}} = \frac{i}{m}+ \frac{j}{m} + \frac{k}{m}=1$ and using the H\"older inequality we get:

$A_{i,j,k} \leq \left\|(g(\varepsilon_{t-1}))^i\right\|_{\frac{m}{i}} \left\|(u(x_{t-1}))^j\right\|_{\frac{m}{j}} \left\|(c(\varepsilon_{t-1}))^i\right\|_{\frac{m}{k}} < \infty$ so that $B < \infty$.

Since $0 \ <  B  < \ \infty $ and $B$ does not contain $x$, we have
\begin{eqnarray*}
E\left[f(\sigma_t^{\delta})|\sigma_{t-1}^{\delta}=x\right] &\leq& 1 + E[c(\varepsilon_{t-1})]^m x^m + (1+x)^{m-1}B= f(x)\left[1- \tau \underbrace{\left(\frac{ x^m - \frac{B}{\tau} (1+x)^{m-1}}{1+x^m}\right)}_{=:C(x)} \right] ,
\end{eqnarray*}
where $\tau:=1-E[c(\varepsilon_{t-1})]^m$.

Since $C(x) \longrightarrow 1$  as $x \longrightarrow +\infty$, there exists a positive real number $K$ such that for $ x \in [K,\infty):  \ \ C(x)>1/2$. Hence for all $x \in (K, \infty):$
\begin{equation}
E\left[f(\sigma_t^{\delta})|\sigma_{t-1}^{\delta}=x\right] \leq (1-\tau/2)f(x) \ .
\label{eq: first help}
\end{equation}

Furthermore for all $x \in (0,K]:$
\begin{equation}
\ E\left[f(\sigma_t^{\delta})|\sigma_{t-1}^{\delta}=x\right] \leq 1+K^m+B(1+K)^{m-1} \ .
\label{eq: second help}
\end{equation}
The equations \eqref{eq: first help} and \eqref{eq: second help} imply that the conditions \eqref{eq:lemma Tweedie 1} and \eqref{eq:lemma Tweedie 2} (with A=(0,K]) of the Tweedie's lemma \ref{lemma of Tweedie} are fulfilled. Thus $\{\sigma_t^{\delta}\}_t$ admits an invariant probability measure $\pi$.

By setting $g(x):=\tau f(x)/2$, we obtain from the inequality \eqref{eq: first help} and the assertion \eqref{eq:lemma Tweedie 3} of the lemma of Tweedie that $E_{\pi}\sigma_t^{m\delta}<\infty$ \ .

Since $\alpha < 1 \leq m$ it follows that $E_{\pi}\sigma_t^{\alpha \delta}<\infty$, i.e., $\pi$ is the $\alpha \delta$th-order stationary distribution of $\sigma_t$ and so $R_t=\sigma_t \varepsilon_t$ is the $\alpha \delta$th-order stationary solution of our model (equations \eqref{eq:first equat of the model} and \eqref{eq:second equat of the model}). Since there is a unique solution of the model (see proposition \ref{propoStatErgoInfty}) we can assert that the solution $R_t$ given in the proposition has a finite moment of order $m$.

This proves the proposition. $ \hspace{10cm} \square$
\newline

\textbf{\textit{Proof of the proposition \ref{geometricergodicity} and \ref{geoErgoExistenceOfV} }}
By considering the noise process as the two dimensional process $(\varepsilon_t,x_t)$, we are in the same situation as in \cite{Kris2005} and can use their proofs for theorem 7 and corollary 8.  $ \hspace{3cm} \square$
\newline


\textbf{\textit{Proof of the lemma \ref{lemma: lemma1 for consistence of QMLE}}}
Due to the fact that the optimum points $\theta_i^*$ belong to $\Theta$ and the fact that all elements of $\Theta$ fulfilled the assumptions of the proposition \ref{propoStatErgoInfty} we have 
\begin{equation*}
E(\sup\limits_{\theta\in\Theta}c(\varepsilon_t; \theta)^{\alpha}) = \max\limits_{i=1,...,n} E(c(\varepsilon_t; \theta_i^*)^{\alpha}) < 1  \ \textrm{ and } \  
E(\sup\limits_{\theta\in\Theta}g(\varepsilon_t; \theta)^{\alpha}) = \max\limits_{i=1,...,n} E(g(\varepsilon_t; \theta_i^*)^{\alpha}) < \infty  \ .
\end{equation*} 
We also have -since $\Theta$ is compact- $\lambda \leq \bar{\lambda} < \infty$, from which follows $E(\sup\limits_{\theta\in\Theta}\lambda u_1(x_{t-1}^{\alpha})) = E(\bar{\lambda} u_1(x_{t-1})^{\alpha}) <\infty$ by assumption $(i)$.

Using the general solution form of the volatility given in the equation \eqref{eq: general form of the solution} and proceeding as in the proof of the proposition \ref{propoStatErgoInfty}, we prove equation \eqref{eq:first equat of lemma1 for consistence of QMLE}.

The equation \eqref{eq:second equat of lemma1 for consistence of QMLE} is proven by the fact that -as stated above-:
\begin{equation*}
0 \leq E(\sup\limits_{\theta\in\Theta}c(\varepsilon_t; \theta)^{\alpha}) = \underbrace{\max\limits_{i=1,...,n} E( c(\varepsilon_t; \theta_i^*)^{\alpha})}_ {=: \ q} < 1 \ .  \hspace{3cm} \square
\end{equation*}


\textbf{\textit{Proof of the lemma \ref{lemma: difference between the both squared volatilities functions}}}
Repeating the proof of the proposition \ref{propoStatErgoTime0} we get that the difference between the observed and the unobserved volatilities vanishes with an exponential rate for every $\theta$. Precisely we get
\begin{equation*}
E\sup\limits_{\theta\in\Theta}\left|\sigma_t^{\delta}(\theta)-\tilde{\sigma}_t^{\delta}(\theta)\right|^{\alpha} = \left( E\sup\limits_{\theta\in\Theta}\underbrace{\left|\tilde{\sigma}_0\right|^{\delta\alpha}}_{A_1(\theta)=\omega^{\alpha \delta}}+ E\sup\limits_{\theta\in\Theta}\underbrace{\left|\sigma_0(\theta)\right|^{\delta \alpha}}_{A_2(\theta)}\right) \left(E\sup\limits_{\theta\in\Theta}\underbrace{\prod\limits_{j=0}^{t-1}[c(\varepsilon_j;\theta)]^{\alpha}}_{A_3(\theta)}\right) \ .
\end{equation*}
It is clear that $A_1$ is continuous in $\theta$. Since $\Theta$ is compact we get $E\left(\sup\limits_{\theta \in \Theta}A_1(\theta)\right)<\infty$.
Using the fact that the function $c$ is non negative, the innovation process $\{\varepsilon_t\}_t$ is $i.i.d.$ and the equation \eqref{eq:second equat of lemma1 for consistence of QMLE} we get  $E\left(\sup\limits_{\theta \in \Theta} A_3(\theta)\right) =  \prod\limits_{j=0}^{t-1}\left[E\sup\limits_{\theta\in\Theta}(c(\varepsilon_j;\theta)^{\alpha})\right] = \rho_1^t$. From the above lemma \ref{lemma: lemma1 for consistence of QMLE}, namely equation \eqref{eq:first equat of lemma1 for consistence of QMLE} we also get that  $E\left(\sup\limits_{\theta \in \Theta}A_2(\theta)\right)<\infty$. Hence we have
\begin{equation*}
E \left(\sup\limits_{\theta \in \Theta} \left|\sigma_t^{\delta}-\tilde{\sigma}_t^{\delta}\right|^{\alpha}\right)= \mathcal{O}(\rho_1^t) \ .
\end{equation*}

To prove the second part of the lemma we first remark that for all $x\in\mathbb{R}$ and $a\in\mathbb{R}_{>0}, \ a\leq1$, we have $\left|x^a-1\right|\leq\left|x-1\right|$. Thus if $2\leq\delta$, i.e., $2/\delta\leq1$ it follows
\begin{eqnarray}
E\sup\limits_{\theta\in\Theta} \left|\sigma_t^{2}(\theta)-\tilde{\sigma}_t^{2}(\theta)\right|^{\alpha} &=& E \left[\sup\limits_{\theta\in\Theta} \tilde{\sigma}_t^{2\alpha}(\theta) \left|\left(\frac{\sigma_t^{\delta}(\theta)}{\tilde{\sigma}_t^{\delta}(\theta)}\right)^{2/\delta}-1\right|^{\alpha}\right] \nonumber \\
&\underbrace{\leq}_{2/\delta \leq1}& E \left[\sup\limits_{\theta\in\Theta} \tilde{\sigma}_t^{2\alpha}(\theta) \left|\frac{\sigma_t^{\delta}(\theta)}{\tilde{\sigma}_t^{\delta}(\theta)}-1\right|^{\alpha}\right] \nonumber \\
&=& E \left[ \sup\limits_{\theta\in\Theta} \underbrace{\tilde{\sigma}_t^{\alpha(2-\delta)}(\theta)}_{\leq \underline{\omega}^{2\alpha(2-\delta)/\delta}=:A_4} \left|\sigma_t^{\delta}(\theta)-\tilde{\sigma}_t^{\delta}(\theta)\right|^{\alpha}\right] \nonumber \\
&=& A_4 E \left[ \sup\limits_{\theta\in\Theta}  \left|\sigma_t^{\delta}(\theta)-\tilde{\sigma}_t^{\delta}(\theta)\right|^{\alpha}\right] \nonumber \\
&=& \mathcal{O}(\rho_1^t) \  \textrm{ by equation } \eqref{eq:difference between observed and unobserved volatilities}.
\end{eqnarray}
If in contrary we have $a>1$, then for all natural number $n+1 \in \mathbb{N}$ greater than $a$ (especially for $n=\left\lfloor a\right\rfloor$) we have $\left|x^a-1\right|\leq\left|x^{n+1}-1\right|$. Thus if $2>\delta$ i.e. $2/\delta>1$ it follows\footnote{For simplicity and a better readability we will omit the dependence on $\theta$.}
\begin{eqnarray}
E\sup\limits_{\theta\in\Theta} \left|\sigma_t^{2}-\tilde{\sigma}_t^{2}\right|^{v} &=& E \left[\sup\limits_{\theta\in\Theta} \tilde{\sigma}_t^{2v} \left|\left(\frac{\sigma_t^{\delta}}{\tilde{\sigma}_t^{\delta}}\right)^{2/\delta}-1\right|^{v}\right] \nonumber \\
&\underbrace{\leq}_{n=\left\lfloor  2/\delta \right\rfloor}& E \left[\sup\limits_{\theta\in\Theta} \tilde{\sigma}_t^{2v} \left|\left(\frac{\sigma_t^{\delta}}{\tilde{\sigma}_t^{\delta}}\right)^{n+1}-1\right|^{v}\right] \nonumber \\
&=& E \left[\sup\limits_{\theta\in\Theta} \underbrace{\left(\tilde{\sigma}_t\right)^{2v-n-1}}_{\leq \underline{\omega}^{2\alpha(n+1-v)/\delta}=:A_5} \left|\left(\sigma_t^{\delta}\right)^{n+1}-\left(\tilde{\sigma}_t^{\delta}\right)^{n+1}\right|^{v}\right] \nonumber \\
&\leq& A_5 \cdot E \left[\sup\limits_{\theta\in\Theta} \left|\sigma_t^{\delta}-\tilde{\sigma}_t^{\delta}\right|^{v} \sum\limits_{j=1}^{n} \left|\sigma_t^{\delta}\left|^{jv} \right|\tilde{\sigma}_t^{\delta}\right|^{(n-j)v}\right] \nonumber \\
&\leq& A_5 \cdot \left[ \underbrace{E\left(\sup\limits_{\theta\in\Theta} \left|\sigma_t^{\delta}-\tilde{\sigma}_t^{\delta}\right|^{v} \right)^2}_{A_6}\right]^{1/2}
 \left[\underbrace{E\left(\sup\limits_{\theta\in\Theta} \sum\limits_{j=1}^{n} \left|\sigma_t^{\delta}\left|^{jv} \right|\tilde{\sigma}_t^{\delta}\right|^{(n-j)v}\right)^2}_{A_7}\right]^{1/2} \ . \nonumber 
\end{eqnarray}

It suffices to choose any $v$ such that $0<v\leq \frac{\alpha}{4n}$ to get the desired result since the following holds: concerning $A_6$ we have $0<v\leq \alpha/2$ and using the equation \eqref{eq:difference between observed and unobserved volatilities} and the Jensen inequality we obtain: $A_6 \leq \left(E \left(\sup\limits_{\theta \in \Theta} \left|\sigma_t^{\delta}-\tilde{\sigma}_t^{\delta}\right|^{\alpha}\right)\right)^{\alpha/2v} = \mathcal{O}((\rho_1^{\alpha/2v})^t)\ \ $. We now prove that $A_{7}$ is finite. Due to the Cauchy Schwartz inequality, sufficient conditions for the finiteness of $A_7$ are that $E \left[\sup\limits_{\theta\in\Theta} \left|\sigma_t^{\delta}\right|^{4jv} \right] < \infty$ and $E\left[\sup\limits_{\theta\in\Theta}\left|\tilde{\sigma}_t^{\delta}\right|^{4(n-j)v}\right] < \infty$. Using the Jensen inequality (since $4jv \  \leq \alpha$, $4(n-j)v \leq \alpha$) together with the fact that $\tilde{\sigma}_t$ is continuous and \eqref{eq:first equat of lemma1 for consistence of QMLE}, we show the two above inequalities. So we set $\rho_2=\rho_1^{\alpha/2v}$ and get the desired result. $\hspace{3cm} \square$
\newline

\textbf{\textit{Proof of the lemma \ref{lemma: difference between the both likelihood functions}}}

Using the value $v$ found in the preceding lemma we have
\begin{eqnarray*}
E\left(n^{v/2}\sup\limits_{\theta \in \Theta}\left|L_n(\theta)-\tilde{L}_n(\theta)\right|^{v/2}\right) &\leq& n^{v/2} \cdot \frac{1}{n^{v/2}}E\left(\sup\limits_{\theta\in\Theta}\sum\limits_{t=1}^{n}\left|\ell_t(\theta)-\tilde{\ell}_t(\theta)\right|^{v/2}\right)  \\
  &=& \sum\limits_{t=1}^{n} E\left(\sup\limits_{\theta\in\Theta}\left|-log(\sigma_t^2)-\frac{R_t^2}{\sigma_t^2}+ log(\tilde{\sigma}_t^2)+\frac{R_t^2}{\tilde{\sigma}_t^2}\right|^{v/2}\right)   \\
  &\leq& \sum\limits_{t=1}^{n} \underbrace{E\left(\sup\limits_{\theta\in\Theta}\frac{\left|\sigma_t^2(\theta)-\tilde{\sigma}_t^2(\theta)\right|^{v/2}}{\sigma_t^{2v/2}\tilde{\sigma}_t^{2v/2}} R_t^{2v/2}\right)}_{B_{1,t}}  \\
&+& \sum\limits_{t=1}^{n} \underbrace{E\left(\sup\limits_{\theta\in\Theta} \left|\log\left(\frac{\sigma_t^2(\theta)} {\tilde{\sigma}_t^2(\theta)}\right)\right|^{v/2} \right)}_{B_{2,t}},
\end{eqnarray*}
noting that the first and the the second inequalities are due to the triangular inequality and the fact that $v/2<1$. Let us now study $B_{1,t}$ and $B_{2,t}$.
\begin{eqnarray}
B_{2,t}&=& E\left(\sup\limits_{\theta\in\Theta} \left|\log\left(\frac{\sigma_t(\theta)} {\tilde{\sigma}_t(\theta)}\right)^2\right|^{v/2} \right) 
= 
E\left(\sup\limits_{\theta\in\Theta} \left|\left(\frac{2}{\delta}\right)\log\left(\frac{\sigma_t(\theta)} {\tilde{\sigma}_t(\theta)}\right)^{\delta}\right|^{v/2} \right) \nonumber \\
&\leq& \left(\frac{2}{\delta}\right)^{v/2} E\left(\sup\limits_{\theta\in\Theta} \left|\frac{\sigma_t^{\delta}(\theta)-\tilde{\sigma}_t^{\delta}(\theta)} {\tilde{\sigma}_t^{\delta}(\theta)}\right|^{v/2}\right) \nonumber \\
&\leq&\left(\frac{2}{\delta}\right)^{v/2} \frac{1}{\underline{\omega}^{v/2}} E\left(\sup\limits_{\theta\in\Theta}\left|\sigma_t^{\delta}(\theta)-\tilde{\sigma}_t^{\delta} (\theta)\right|^{v/2} \right) \nonumber \\
&\leq&\left(\frac{2}{\delta}\right)^{v/2} \frac{1}{\underline{\omega}^{v/2}} \left[\underbrace{E\left(\sup\limits_{\theta\in\Theta}\left|\sigma_t^{\delta}(\theta)-\tilde{\sigma}_t^{\delta} (\theta)\right|^{v/2} \right)}_{= \ \mathcal{O}(\rho_1^t) \textrm{ by equation  \eqref{eq:difference between observed and unobserved volatilities}} }\right]^{v/\alpha} = \mathcal{O}\left(\left[\rho_1^{v/\alpha}\right]^t\right),\nonumber 
\end{eqnarray}
noting that the first inequality is due to the fact that $\log(x)\leq x-1$ and the last one follows from the Jensen inequality using the fact that $v/2 \leq \alpha$. Note also that $0<\rho_1^{v/\alpha}<1$, since $0<\rho_1<1$. 

Concerning $B_{1,t}$ we use the equation \eqref{lemma: difference between the both squared volatilities functions} in the following sense
\begin{eqnarray}
B_{1,t}&=&E\left(\sup\limits_{\theta\in\Theta}\frac{\left|\sigma_t^2(\theta)-\tilde{\sigma}_t^2(\theta)\right|^{v/2}}{\sigma_t^{v}(\theta)\tilde{\sigma}_t^{v}(\theta)} R_t^{v}\right) \nonumber \\ 
&\leq&\left[E\left(\sup\limits_{\theta\in\Theta}\frac{\left|\sigma_t^2(\theta)-\tilde{\sigma}_t^2(\theta)\right|^{v}}{\sigma_t^{2v}(\theta)\tilde{\sigma}_t^{2v}(\theta)}\right)\right]^{1/2} \cdot \left[E\left(\sigma_t^{2v}(\theta_0)\varepsilon_t^{2v}\right)\right]^{1/2} \nonumber \\
&=&\left[E\left(\sup\limits_{\theta\in\Theta}\frac{\left|\sigma_t^2(\theta)-\tilde{\sigma}_t^2(\theta)\right|^{v}}{\sigma_t^{2v}(\theta)\tilde{\sigma}_t^{2v}(\theta)}\right)\right]^{1/2} \cdot \left[E\sigma_t^{2v}(\theta_0)\right]^{1/2} E\left[\varepsilon_t^{2v}\right]^{1/2} \nonumber \\
&\leq& \frac{1}{\underline{\omega}^{4v/\delta}}\left[\underbrace{E\left(\sup\limits_{\theta\in\Theta}\left|\sigma_t^2(\theta)-\tilde{\sigma}_t^2(\theta)\right|^{v}\right)}_{= \ \mathcal{O}(\rho_2^t) \textrm{ by equation  \eqref{eq:difference between squared observed and unobserved volatilities}}} \right]^{1/2} \underbrace{\left[E\sigma_t^{2v}(\theta_0)\right]^{1/2} E\left[\varepsilon_t^{2v}\right]^{1/2}}_{<\infty \ \textrm{by assumption \ref{ass: assumption for consistency} (ii)} } =\mathcal{O}\left(\left[\rho_2^{1/2}\right]^t\right). \nonumber 
\end{eqnarray}
Note that the first inequality relies on the Cauchy Schwartz inequality and the second equality (decomposition of the expectation) is due to the fact the $\varepsilon_t$ is independent of $\sigma_t(\theta)$.

We get:
\begin{eqnarray*}
E\left(n^{v/2}\sup\limits_{\theta \in \Theta}\left|L_n(\theta)-\tilde{L}_n(\theta)\right|^{v/2}\right) &\leq& \sum\limits_{t=1}^{n} \mathcal{O}\left(\left[\rho_2^{1/2}\right]^t\right) + \mathcal{O}\left(\left[\rho_1^{v/\alpha}\right]^t\right) \\
&\leq& \sum\limits_{t=1}^{\infty}\mathcal{O}\left(\left[\rho_2^{1/2}\right]^t\right)+\mathcal{O}\left(\left[\rho_1^{v/\alpha}\right]^t\right)=:C < \infty 
\end{eqnarray*}
and $C$ is independent of n.

Finally we get (through the Markov inequality) for every positive $\epsilon>0$:
\begin{eqnarray*}
\mathbb{P}\left(\sup\limits_{\theta \in \Theta}\left|L_n(\theta)-\tilde{L}_n(\theta)\right| > \epsilon\right) &=& \mathbb{P}\left(n^{v/2}\sup\limits_{\theta \in \Theta}\left|L_n(\theta)-\tilde{L}_n(\theta)\right|^{v/2} > n^{v/2} \epsilon^{v/2}\right) \\
&\leq& \frac{E\left(n^{v/2}\sup\limits_{\theta \in \Theta}\left|L_n(\theta)-\tilde{L}_n(\theta)\right|^{v/2}\right)}{n^{v/2} \epsilon^{v/2}} \\ 
&\leq& \frac{C}{{n^{v/2} \epsilon^{v/2}}} \ \ \longrightarrow 0 \ \textrm{ as } n\rightarrow\infty  \ . \hspace{3cm} \square
\end{eqnarray*}

\textbf{\textit{Proof of the lemma \ref{lemma: lemma0 for consistence}}}

To prove this lemma we will mainly use the last assumption $(v)$ together with the fact that $\{x_t\}_t$ is independent of $\{\varepsilon_t\}_t$. These imply -since u is non-constant- that $x_t|\left\{\mathfrak{F}_{t-1}, \{\varepsilon_t\}_t \right\}$ has a non-degenerate distribution.  

We have  $\sigma_t^{\delta}(\theta)= \sum\limits_{k=0}^{\infty} A_k(\theta)$ where $A_k(\theta)=\left[g(\varepsilon_{t-1-k};\ \omega, \theta^1)+\lambda u_1(x_{t-1-k})\right] \prod\limits_{j=0}^{k-1} c(\varepsilon_{t-1-j}; \theta^2) \ .$
\begin{eqnarray*}
\sigma_t^{\delta}(\theta)=\sigma_t^{\delta}(\theta_0) & \Rightarrow & A_0(\theta)-A_0(\theta_0) = \sum\limits_{k=0}^{\infty} \left(A_k(\theta_0)-A_k(\theta)\right) \\
& \Rightarrow & (\lambda - \lambda_0)u(x_{t-1})=g(\varepsilon_{t-1-k};\ \omega_0, \theta_0^1)-g(\varepsilon_{t-1-k};\ \omega, \theta^1) + \sum\limits_{k=0}^{\infty}\left(A_k(\theta_0)-A_k(\theta)\right)  \ .
\end{eqnarray*}

Since $x_t|\left\{\mathfrak{F}_{t-1}, \{\varepsilon_t\}_t \right\}$ has a non-degenerate distribution we get $\lambda = \lambda_0$. Applying the same argument for $A_k(\theta)-A_k(\theta_0)$, \ $k=1,2,3,...$ we get that $\prod\limits_{j=0}^{k-1} c(\varepsilon_{t-1-j}; \theta_0^2)=\prod\limits_{j=0}^{k-1} c(\varepsilon_{t-1-j}; \theta^2)$ for $k=1,2,3,...$, especially $c(\varepsilon_{t-1-j}; \theta_0^2) =c(\varepsilon_{t-1-j}; \theta^2)$.  Finally plugging this in the equality $\sigma_t^{\delta}(\theta)=\sigma_t^{\delta}(\theta_0)$ we get $\sum\limits_{k=0}^{\infty} \left[g(\varepsilon_{t-1-k};\ \omega, \theta^1)- g(\varepsilon_{t-1-k};\ \omega_0, \theta_0^1)\right] \prod\limits_{j=0}^{k-1} c(\varepsilon_{t-1-j}; \theta_0^2)=0$
Since $c(\varepsilon_{t-1-j}; \theta_0^2)>0$ we get $g(\varepsilon_{t-1-k};\ \omega, \theta^1)- g(\varepsilon_{t-1-k};\ \omega_0, \theta_0^1)$.

From the equalities $c(\varepsilon_{t-1-j}; \theta^2)=c(\varepsilon_{t-1-j}; \theta_0^2)$, $g(\varepsilon_{t-1-k};\ \omega, \theta^1)- g(\varepsilon_{t-1-k};\ \omega_0, \theta_0^1)$ and the assumption $(iii)$, we get the desired result.
$\hspace{3cm} \square$
\newline

\textbf{\textit{Proof of the lemma \ref{lemma: lemma 2 for the proof of consistency}}}

Since $\frac{R_t^2}{\sigma_t^2(\theta_0)}>0$, we have on the one hand:
\begin{eqnarray*}
\ell_t(\theta_0)= -log(\sigma_t^2(\theta_0))-\frac{R_t^2}{\sigma_t^2(\theta_0)} \leq - log(\sigma_t^2(\theta_0)) \leq -log(\omega_0) \ ,
\end{eqnarray*}
so that $E\left[\left(\ell_t(\theta_0)\right)^{+}\right] < \infty \ .$

On the other hand we have
\begin{eqnarray*}
\ell_t(\theta_0)=-log(\sigma_t^2(\theta_0))-\frac{\sigma_t^2(\theta_0) \varepsilon_t^2 }{\sigma_t^2(\theta_0)} = -log(\sigma_t^2(\theta_0))-\frac{\sigma_t^2(\theta_0) \varepsilon_t^2 }{\sigma_t^2(\theta_0)}
\end{eqnarray*}
and since $E(\varepsilon_t^2)<\infty$ (by assumption \ref{ass: assumption for consistency} 2.), we obtained through the Jensen inequality and the first assertion of assumption \ref{ass: assumption for consistency} that the negative part is bounded.
\begin{eqnarray*}
E\left[\left(\log\sigma_t^2(\theta_0)\right)^{-}\right]= \frac{1}{s}E\left[\left(\log\sigma_t^{2s}(\theta_0)\right)^{-}\right] \leq \frac{1}{s}\left(E\left[\log\sigma_t^{2s}(\theta_0)\right]\right)^{-} < +\infty \ ,
\end{eqnarray*}
so that the expectation exists. This concludes the first part of the lemma. For the second part of the lemma we first remark that as proved above $E[\ell_t(\theta)]^{+}< \infty$, such that $E[\ell_t(\theta)] \in [-\infty, \infty)$.

If $E[\ell_t(\theta)]= -\infty$ then $E[\ell_t(\theta)]<E[\ell_t(\theta_0)] \ .$

If $E[\ell_t(\theta)] > -\infty$ we have
\begin{eqnarray}
E[\ell_t(\theta_0)]-E[\ell_t(\theta)] &=& E\left[-log(\sigma_t^2(\theta_0))-\frac{R_t^2}{\sigma_t^2(\theta_0)}+ log(\sigma_t^2(\theta))+\frac{R_t^2}{\sigma_t^2(\theta)}\right] \nonumber \\
&=&E\left[\log\left(\frac{\sigma_t^2(\theta)}{\sigma_t^2(\theta_0)}\right)+ \frac{\sigma_t^2(\theta_0)}{\sigma_t^2(\theta)} \right] -1 \geq 0 \ ,
\end{eqnarray}
since for all $x >0, \  \ x+\log(1/x) \geq 1$.

The equality occurs if and only if $\sigma_t^2(\theta_0)=\sigma_t^2(\theta)$ or $\sigma_t^{\delta}(\theta_0)=\sigma_t^{\delta}(\theta)$ which implies from the previous lemma that $\theta=\theta_0.$
$ \hspace{3cm} \square$
\newline


\textbf{\textit{Proof of the lemma \ref{lemma: lemma 3 for the proof of consistency}}}

Let $\theta \in \Theta \setminus \{\theta_0\}$ and denote by $V_{1/k}(\theta)$ the open centered in $\theta$ with radius $1/k$. 
Since $\{\ell_t(\theta)\}_t$ is strictly stationary and ergodic and $E{\ell_t(\theta)}\in \mathbb{R} \cup \{-\infty\}$ we can use a modified  version of the ergodic theorem\footnote{see \cite{FranZako2010}, Exercises 7.3 and 7.4.} on $\{\ell_t(\theta)\}_t$ and thus on $\left\{\sup\limits_{\theta^* \in V_{1/k}(\theta)\cap \Theta} \ell_t(\theta^*)\right\}_t$ and we get on one side:
\begin{eqnarray}
\limsup\limits_{n\rightarrow \infty} \sup\limits_{\theta^* \in V_{1/k}(\theta)\cap \Theta} L_n(\theta^*)=E_{\theta_0} \left(\sup\limits_{\theta^* \in V_{1/k}(\theta)\cap \Theta} \ell_t(\theta^*)\right) <  E_{\theta_0}\left(\ell_t(\theta_0)\right) ,
\label{eq: first inequality}
\end{eqnarray}
where the inequality is obtained by the following argument: using Beppo-Levi's theorem, the fact that $\left(\sup\limits_{\theta^* \in V_{1/k}(\theta)\cap \Theta} \ell_t(\theta^*)\right)_{k\in\mathbb{N}}$ is a decreasing sequence converging through $\ell_t(\theta)$ and since $E_{\theta_0}\left(\ell_t(\theta)\right) <  E_{\theta_0}\left(\ell_t(\theta_0)\right)$, there exists a $k \in \mathbb{N}$ such that the above inequality holds. 

On the other side we have for any neighborhood $V(\theta_0)$ of $\theta_0$
\begin{eqnarray}
\liminf\limits_{n\rightarrow \infty} \sup\limits_{\theta^* \in V(\theta_0)\cap \Theta} L_n(\theta^*)=E_{\theta_0} \left(\sup\limits_{\theta^* \in V(\theta_0)\cap \Theta} \ell_t(\theta^*)\right) \geq E_{\theta_0}\left(\ell_t(\theta_0)\right),
\label{eq: second inequality}
\end{eqnarray}
where the inequality is due to the fact that $\theta_0 \in V(\theta_0)$.

Let $k\in \mathbb{N}$ and $V(\theta_0)$ be a neighborhood of $\theta_0$. Since $\Theta$ is compact, there exists a finite number of elements $\theta_1,...,\theta_{n_k}$ different from $\theta$  such that $\bigcup\limits_{i=1}^{n_k}V(\theta_i) \!\supseteq \Theta$, where $V(\theta_i), \ i=1,...,k$ is an open ball with radius $1/k$. Applying the inequalities \eqref{eq: first inequality} and \eqref{eq: second inequality}, we get
\begin{eqnarray}
\lim\limits_{n\rightarrow \infty}  \sup\limits_{\theta^* \in \Theta} L_n(\theta^*) &=& \lim\limits_{n\rightarrow \infty} \  \sup\limits_{i=0,...,n_k} \sup\limits_{\theta^* \in V(\theta_i) \cap \Theta} \ L_n(\theta^*) \nonumber \\
&=& E_{\theta_0}\left(\ell_t(\theta_0)\right) \ \ a.s.
\label{eq: main equation of lemma 3}
\end{eqnarray}
This proves the lemma. $ \hspace{10cm} \square$
\newline

We can now show the consistency of the Q-MLE.


\textbf{\textit{Proof of Proposition \ref{prop: consistency of the QMLE}}}

The proof of the proposition \ref{prop: consistency of the QMLE} relies on the Theorem 4.1.1 of \cite{Amen1985}.

From lemmas \ref{lemma: difference between the both likelihood functions} and \ref{lemma: lemma 3 for the proof of consistency} we get
\begin{eqnarray*}
\sup\limits_{\theta\in\Theta}\left|\tilde{L}_n(\theta)-E_{\theta_0}\left(\ell_t(\theta_0)\right)\right| &\leq& \sup\limits_{\theta\in\Theta}\left|\tilde{L}_n(\theta)-L_n(\theta)\right| +\sup\limits_{\theta\in\Theta}\left|L_n(\theta)-E_{\theta_0}\left(\ell_t(\theta_0)\right)\right| \longrightarrow \ 0 \ (\textrm{in probability}).
\end{eqnarray*}
Furthermore lemma \ref{lemma: lemma 2 for the proof of consistency} shows that the function $\theta\ \longrightarrow \ E_{\theta_0}\left(\ell_t(\theta)\right)$ admits a global maximum at $\theta_0.$  Hence the proposition is proven.$ \hspace{10cm} \square$
\newline

 
\textbf{\textit{Proof of lemma \ref{lemma: first lemma for ass norm}}}

Since
\begin{equation}
\frac{\partial \ell_t}{\partial \theta }(\theta)=\frac{2}{\delta}\frac{1}{\sigma_t^{\delta}(\theta)} \frac{\partial \sigma_t^{\delta}(\theta)}{\partial \theta }(\theta) \left[\frac{R_t^2}{\sigma_t^2(\theta)}-1\right],
\label{eq: def of deriv of ell}
\end{equation}
we have
\begin{eqnarray*}
E_{\theta_0}\left(\frac{\partial \ell_t}{\partial \theta }(\theta_0)\right) &=&E_{\theta_0}\left(\frac{2}{\delta}\frac{1}{\sigma_t^{\delta}(\theta_0)} \frac{\partial \sigma_t^{\delta}(\theta_0)}{\partial \theta } (\theta_0) \left[\frac{\sigma_t^2(\theta_0)\varepsilon_t^2}{\sigma_t^2(\theta_0)}-1\right]\right) \\
&=&
E_{\theta_0}\left(\frac{2}{\delta}\frac{1}{\sigma_t^{\delta}(\theta_0)} \frac{\partial \sigma_t^{\delta}(\theta_0)}{\partial \theta } (\theta_0) \right) \underbrace{E_{\theta_0}\left(\varepsilon_t^2 - 1\right)}_{=0} \  \ =  \ 0 \ .
\end{eqnarray*}
The second equality is due to the fact that $\varepsilon_t$ is independent of $\mathfrak{F}_{t-1}$.

Since the expectation is null the variance equals the second moment; we then have
\begin{eqnarray*}
var\left(\frac{\partial \ell_t}{\partial \theta }(\theta_0)\right)&=&E_{\theta_0}\left(\frac{\partial \ell_t(\theta_0)}{\partial\theta } \frac{\partial \ell_t(\theta_0)}{\partial \theta^{'} }\right) \\
&=& E_{\theta_0}\left(\frac{4}{\delta^2} \frac{1}{\sigma_t^{2\delta}}\frac{\partial \sigma_t^{\delta}(\theta_0)}{\partial\theta }  \frac{\partial \sigma_t^{\delta}(\theta_0)}{\partial \theta^{'} }\right) E_{\theta_0}\left(\left(\varepsilon_t^2 - 1\right)^2\right) \\
&=&  \frac{4}{\delta^2}E_{\theta_0}\left( \frac{1}{\sigma_t^{2\delta}}\frac{\partial \sigma_t^{\delta}(\theta_0)}{\partial\theta }  \frac{\partial \sigma_t^{\delta}(\theta_0)}{\partial \theta^{'} }\right) E_{\theta_0}\left(\varepsilon_t^4 -2 \varepsilon_t^2 + 1^2\right) \\
&=&\frac{4}{\delta^2}\left(\kappa-1\right) A \ .  \hspace{3cm} \square
\end{eqnarray*}

 
\textbf{\textit{Proof of the lemma \ref{lemma: second lemma for ass norm}}}

In equation \eqref{eq: first derivative of sigma}, $E(M_1^{\alpha})<\infty$ by assumption \ref{ass: assumptions for the ass norm} $(iv)$. $E(M_2^{\alpha})<\infty$ since $E(u(x_t)^{\alpha})<\infty$ as it can be read in proposition \ref{propoStatErgoInfty}. By Cauchy Schwartz inequality we have $E(M_3^{\alpha/2})=E\left|\frac{\partial c}{\partial\theta_i}\sigma_{t-1}^{\delta} \right|^{\alpha/2} \leq  \left(E\left|\frac{\partial c}{\partial\theta_i}\right|^{\alpha} E\sigma_{t-1}^{\delta\alpha} \right)^{1/2} <\infty$. Thus using the same procedure as in the proof of proposition \ref{propoStatErgoInfty}, we obtain that $\frac{\partial\sigma_t^{\delta}}{\partial \theta}$ is $\alpha/2$ stationary and ergodic. Similarly it can be shown that $\frac{\partial^2 \sigma_t^{\delta}} {\partial\theta\partial \theta^{'}}$ is a $\alpha/4$ stationary and ergodic process. We note that the stationarity and ergodic properties stated above follow from the fact that $\partial\sigma_t^{\delta}/\partial\theta$ and $\partial^2 \sigma_t^{\delta}/\partial\theta\partial\theta^{'}$ are fixed functionals of the joint process $\{(\varepsilon_t,x_t)\}_t$, which is by assumption \ref{assumption on noise processes} ergodic and stationary.

Looking now at the formulas \eqref{eq: first derivative of ell} and \eqref{eq: second derivative of ell} defining the first and second partial derivatives of $\ell_t$, we can state that they are fixed functional of  $\{(\varepsilon_t,x_t)\}_t$ and hence also stationary and ergodic. $\hspace{3cm} \square$
\newline

 
\textbf{\textit{Proof of the lemma \ref{lemma: third lemma for ass norm}}}
The proof of this lemma is exactly analogous to the proof of the lemma \ref{lemma: difference between the both squared volatilities functions} together with the fact that the volatility $\sigma_t^{\delta}$ and its derivatives are bounded in $V(\theta_0)$ and the fact that $\sigma_t^{\delta}$ (resp. $\partial \sigma_t^{\delta}(\theta)/\partial\theta$, resp. $\partial^2 \sigma_t^{\delta}(\theta)\partial\theta\partial\theta^{'}$) is a $\alpha$ (resp. $\alpha/2$, resp. $\alpha/4$) stationary process.  $\hspace{7cm} \square$
\newline

 
\textbf{\textit{Proof of the lemma \ref{lemma: fourth lemma for ass norm diff ell and tilde ell}}} \newline
To prove \eqref{eq:diff Sn and tilde Sn} we use the form of $\partial\ell_t/\partial\theta$ and $\partial\tilde{\ell}_t/\partial\theta$ given by equation \eqref{eq: first derivative of ell} together with the fact that for $a_1,a_2,a_3,b_1b_2,b_3$ real numbers we have $a_1a_2a_3-b_1b_2b_3=a_2a_3(a_1-b_1)+b_1a_3(a_2-b_2)+b_1b_2(a_3-b_3)$. We get
\begin{eqnarray*}
\sqrt{n}\left(S_n(\theta_0)-\tilde{S}_n(\theta_0)\right)&=& \frac{2}{\sqrt{n}\delta}\sum\limits_{t=1}^{n}\underbrace{\left(\frac{1}{\sigma_t^{\delta}(\theta_0)}-\frac{1}{\tilde{\sigma}_t^{\delta}(\theta_0)}\right)\frac{\partial \sigma_t^{\delta}(\theta_0)}{\partial\theta}\left[\frac{R_t^2}{\sigma_t^{2}(\theta_0)}-1\right]}_{=:D_{1,t}} \\
 && + \frac{2}{\sqrt{n}\delta} \sum\limits_{t=1}^{n} \underbrace{\left(\frac{\partial \sigma_t^{\delta}(\theta_0)}{\partial\theta}- \frac{\partial \tilde{\sigma}_t^{\delta}(\theta_0)}{\partial\theta}\right)\frac{1}{\tilde{\sigma}_t^{\delta}(\theta_0)}\left[\frac{R_t^2}{\sigma_t^{2}(\theta_0)}-1\right]}_{=:D_{2,t}}  \\
 && + \frac{2}{\sqrt{n}\delta} \sum\limits_{t=1}^{n} \underbrace{\frac{1}{\tilde{\sigma}_t^{\delta}(\theta_0)}\frac{\partial \tilde{\sigma}_t^{\delta}(\theta_0)} {\partial\theta} \left[\frac{R_t^2}{\sigma_t^{2}(\theta_0)}-\frac{R_t^2}{\tilde{\sigma}_t^{2}(\theta_0)}\right] }_{=:D_{3,t}} \ . 
\end{eqnarray*}
Let us now study $D_{1,t}$, $D_{2,t}$ and $D_{3,t}$ separately. 

For $D_{1,t}$ we use the fact that $\frac{\partial \sigma_t^{\delta}(\theta_0)}{\partial\theta}$ is $\alpha/2$-stationary, $R_t^2$ is $\alpha$-stationary, $\varepsilon_t$ is independent of $\sigma_t^{\delta}$ and $\tilde{\sigma}_t^{\delta}$, the equation \eqref{eq:difference between observed and unobserved volatilities}, the Jensen and Cauchy-Schwartz inequality to obtain
\begin{eqnarray*}
E_{\theta_0}\left|D_{1,t}\right|^{\alpha/4} &=& E_{\theta_0}\left| \frac{\tilde{\sigma}_t^{\delta}(\theta_0)-\sigma_t^{\delta}(\theta_0)}{\sigma_t^{\delta}(\theta_0)\tilde{\sigma}_t^{\delta}(\theta_0)} \frac{\partial \sigma_t^{\delta}(\theta_0)}{\partial\theta}\left(\varepsilon_t^2-1\right) \right|^{\alpha/4} \\
&=& E_{\theta_0}\left| \frac{\tilde{\sigma}_t^{\delta}(\theta_0)-\sigma_t^{\delta}(\theta_0)}{\underbrace{\sigma_t^{\delta}(\theta_0)}_{\geq \omega_0}\underbrace{\tilde{\sigma}_t^{\delta}(\theta_0)}_{\geq \omega_0}} \frac{\partial \sigma_t^{\delta}(\theta_0)}{\partial\theta}\right|^{\alpha/4} E_{\theta_0}\left|\varepsilon_t^2-1\right|^{\alpha/4}\\
&\leq& \frac{1}{\omega_0^{\alpha/2}} \left( \underbrace{E_{\theta_0}\left| \tilde{\sigma}_t^{\delta}(\theta_0)-\sigma_t^{\delta}(\theta_0)\right|^{\alpha/2}}_{\leq \mathcal{O}\left(\left[\rho_1^{1/2}\right]^t\right)  \textrm{ by equation \eqref{eq:difference between observed and unobserved volatilities}}}\right)^{1/2} \left(\underbrace{E_{\theta_0}\left|\frac{\partial \sigma_t^{\delta}(\theta_0)}{\partial\theta}\right|^{\alpha/2}}_{<\infty}\right)^{1/2} \underbrace{E_{\theta_0} \left|\varepsilon_t^2-1\right|^{\alpha/4}}_{<\infty \textrm{ by assump. \ref{ass: assumptions for the ass norm},2}}\\
& \leq & \mathcal{O}\left(\left[\rho_1^{1/4}\right]^t\right) \ .
\end{eqnarray*}

For $D_{2,t}$ we use almost all arguments listed for $D_{1,t}$ plus the equation \eqref{eq: diff observ and unobserv first der volatility} and we obtain
\begin{eqnarray*}
E_{\theta_0} \left|D_{2,t}\right|^{\alpha/2} &=& E_{\theta_0} \left|\left(\frac{\partial \sigma_t^{\delta}(\theta_0)}{\partial\theta}- \frac{\partial \tilde{\sigma}_t^{\delta}(\theta_0)}{\partial\theta}\right)\frac{1}{\underbrace{\tilde{\sigma}_t^{\delta}(\theta_0)}_{\geq \omega_0}}\left(\varepsilon_t^2-1\right)\right|^{\alpha/2} \\
&\leq&  \frac{1}{\omega_0^{\alpha/2}}\underbrace{E_{\theta_0} \left|\frac{\partial \sigma_t^{\delta}(\theta_0)}{\partial\theta}- \frac{\partial \tilde{\sigma}_t^{\delta}(\theta_0)}{\partial\theta} \right|^{\alpha/2}}_{=\mathcal{O}(\rho_3^t) \textrm{ by equation \eqref{eq: diff observ and unobserv first der volatility}} } \cdot \underbrace{E_{\theta_0} \left|\varepsilon_t^2-1\right|^{\alpha/2}}_{<\infty \textrm{ by assump. \ref{ass: assumptions for the ass norm},2}} \\
& \leq &  \mathcal{O}(\rho_3^t) \ .
\end{eqnarray*}

$D_{3,t}$ is also treated using the above listed arguments together with the equation \eqref{eq: diff observ and unobserv square first der volatility} and we get
\begin{eqnarray*}
E_{\theta_0} \left|D_{3,t}\right|^{v/4} &=&  E_{\theta_0} \left| \frac{1}{\underbrace{\tilde{\sigma}_t^{\delta}(\theta_0)}_{\geq \omega_0}}\frac{\partial \tilde{\sigma}_t^{\delta}(\theta_0)}{\partial\theta}  \left(\frac{\tilde{\sigma}_t^{2}(\theta_0)-\sigma_t^{2}(\theta_0)}{\underbrace{\sigma_t^{2}(\theta_0)}_{\geq \omega_0^{2/\delta}} \underbrace{\tilde{\sigma}_t^{2}(\theta_0)}_{\geq \omega_0^{2/\delta}}}\right) R_t^2 \right|^{v/4} \\
&\leq & \frac{1}{\omega^u} \left(E_{\theta_0} \left|\frac{\partial \tilde{\sigma}_t^{\delta}(\theta_0)}{\partial\theta}\right|^{v/2} \right)^{1/2} \left(E_{\theta_0} \left| \left(\tilde{\sigma}_t^{2}(\theta_0)-\sigma_t^{2}(\theta_0)\right) R_t^2 \right|^{v/2} \right)^{1/2} \\
&\leq & \frac{1}{\omega^u} \left(\underbrace{E_{\theta_0} \left|\frac{\partial \tilde{\sigma}_t^{\delta}(\theta_0)}{\partial\theta}\right|^{v/2} }_{=: D_{3,1,t} < \infty}\right)^{1/2} \left(\underbrace{E_{\theta_0} \left| \tilde{\sigma}_t^{2}(\theta_0)-\sigma_t^{2}(\theta_0) \right|^v}_{= \mathcal{O}(\rho_5^t) \textrm{ by equation \eqref{eq: diff observ and unobserv square first der volatility} }} \right)^{1/4}  \left(\underbrace{E_{\theta_0}R_t^{2v}}_{=: D_{3,2,t} < \infty}\right)^{1/4} \\
&=& \mathcal{O}\left(\left[\rho_5^{1/4}\right]^t\right) ,
\end{eqnarray*}
where $u=v(1/2+2/\delta)$. We note that $D_{3,1,t} < \infty$ is due to the fact that the first derivative of $\sigma_t^{\delta}$ is $\alpha/2$ stationary - as shown in the proof of lemma \ref{lemma: second lemma for ass norm}- and $v/2 < \alpha/2$. To prove $D_{3,2,t} < \infty$ we use the fact that $R_t$ is $\alpha$-stationary and $2v < \alpha$.

Since $\min\{\alpha/2, \alpha/4, v/4\}=v/4$ we get by using Jensen inequality and the above calculations:
\begin{eqnarray*}
E_{\theta_0}\left(n^{v/8}\left|\sqrt{n}\left(S_n(\theta_0)-\tilde{S}_n(\theta_0)\right)\right|^{v/4}\right) &\leq& \left(\frac{2}{\delta}\right)^{v/4}\sum\limits_{t=1}^{n} \left[E_{\theta_0}D_{1,t}^{v/4} + E_{\theta_0}D_{2,t}^{v/4} + E_{\theta_0}D_{3,t}^{v/4}\right]  \\
&\leq& \left(\frac{2}{\delta}\right)^{v/4}\sum\limits_{t=1}^{n} \left[\left(E_{\theta_0}D_{1,t}^{\alpha/4}\right)^{v/\alpha} + \left(E_{\theta_0}D_{2,t}^{\alpha/2}\right)^{v/2\alpha} + E_{\theta_0}D_{3,t}^{v/4}\right]  \\
&\leq&  \left(\frac{2}{\delta}\right)^{v/4}\sum\limits_{t=1}^{\infty} \left[\mathcal{O}\left(\left[\rho_1^{v/4\alpha}\right]^t\right) +  \mathcal{O}\left(\left[\rho_3^{v/2\alpha}\right]^t\right) +  \mathcal{O}\left(\left[\rho_5^{1/4}\right]^t\right)  \right] \\
& =: & F_1 < \infty \ .
\end{eqnarray*}
Thus using the Markov inequality it follows for all $\varepsilon > 0$
\begin{eqnarray*}
\mathbb{P}\left(\sqrt{n}\left|S_n(\theta_0)-\tilde{S}_n(\theta_0)\right| > \varepsilon \right) &=& \mathbb{P}\left(n^{v/8}\left|\sqrt{n}\left(S_n(\theta_0)-\tilde{S}_n(\theta_0)\right)\right|^{v/4} > n^{v/8}\varepsilon^{v/4} \right) \\
& \leq & \frac{E_{\theta_0}\left(n^{v/8}\left|\sqrt{n}\left(S_n(\theta_0)-\tilde{S}_n(\theta_0)\right)\right|^{v/4}\right)}{n^{v/8}\varepsilon^{v/4}} \\
&\leq& \frac{F_1}{n^{v/8}\varepsilon^{v/4}} \  \  \longrightarrow \  0 \  (\textrm{as} \ n \rightarrow \infty).
\end{eqnarray*}


To prove \eqref{eq:diff Hn and tilde Hn} we use the form of $\partial^2\ell_t/\partial\theta\partial\theta^{'}$ and $\partial^2\tilde{\ell}_t/\partial\theta\partial\theta^{'}$ given by equation \eqref{eq: first derivative of ell} together with the fact that for $a_1,a_2,a_3,a_4,b_1,b_2,b_3,b_4$ real numbers we have $a_1a_2a_3a_4-b_1b_2b_3b_4=a_2a_3a_4(a_1-b_1)+b_1a_3a_4(a_2-b_2)+b_1b_2a_4(a_3-b_3)+ b_1b_2b_3(a_4-b_4)$ and $a_1a_2a_3-b_1b_2b_3=a_2a_3(a_1-b_1)+b_1a_3(a_2-b_2)+b_1b_2(a_3-b_3)$. We get for $i,j \in \{1,...,m\}$:
\begin{eqnarray*}
\left(H_n(\theta)-\tilde{H}_n(\theta)\right)_{ij} &=&  \frac{\partial^2 \ell_t}{\partial\theta_i\partial\theta_j} - \frac{\partial^2 \tilde{\ell}_t}{\partial\theta_i\partial\theta_j} \\
&=&\frac{-2}{n\delta}\sum\limits_{t=1}^{n}\underbrace{\left(\frac{1}{\sigma_t^{2\delta}(\theta)}-\frac{1}{\tilde{\sigma}_t^{2\delta}(\theta)}\right)\frac{\partial \sigma_t^{\delta}(\theta)}{\partial\theta_i}\frac{\partial \sigma_t^{\delta}(\theta_i)}{\partial\theta_j} \left[\frac{2+\delta}{\delta} \frac{R_t^2}{\sigma_t^{2}(\theta)}-1\right]}_{=:G_{1,t}^{ij}(\theta)} \\
 && - \frac{2}{n\delta} \sum\limits_{t=1}^{n} \underbrace{\left(\frac{\partial \sigma_t^{\delta}(\theta)}{\partial\theta_i}- \frac{\partial \tilde{\sigma}_t^{\delta}(\theta)}{\partial\theta_i}\right) \frac{1}{\tilde{\sigma}_t^{2\delta}(\theta)}  \frac{\partial \sigma_t^{\delta}(\theta)}{\partial\theta_j} \left[\frac{2+\delta}{\delta} \frac{R_t^2}{\sigma_t^{2}(\theta)}-1\right]}_{=:G_{2,t}^{ij}(\theta)}  \\
 && - \frac{2}{n\delta} \sum\limits_{t=1}^{n} \underbrace{\left(\frac{\partial \sigma_t^{\delta}(\theta)}{\partial\theta_j}- \frac{\partial \tilde{\sigma}_t^{\delta}(\theta)}{\partial\theta_j}\right) \frac{1}{\tilde{\sigma}_t^{2\delta}(\theta)}  \frac{\partial \tilde{\sigma}_t^{\delta}(\theta)}{\partial\theta_i} \left[\frac{2+\delta}{\delta} \frac{R_t^2}{\sigma_t^{2}(\theta)}-1\right]}_{=:G_{3,t}^{ij}(\theta)}  \\
&& - \frac{2(2+\delta)}{n\delta^2} \sum\limits_{t=1}^{n} \underbrace{\left(\frac{R_t^2}{\sigma_t^{2}(\theta)}-\frac{R_t^2}{\tilde{\sigma}_t^{2}(\theta)}\right)\frac{1}{\tilde{\sigma}_t^{2\delta}(\theta)}  \frac{\partial \tilde{\sigma}_t^{\delta}(\theta)}{\partial\theta_i} \frac{\partial \tilde{\sigma}_t^{\delta}(\theta)}{\partial\theta_j} }_{=:G_{4,t}^{ij}(\theta)} \\
&+&\frac{2}{n\delta}\sum\limits_{t=1}^{n}\underbrace{\left(\frac{1}{\sigma_t^{\delta}(\theta)}-\frac{1}{\tilde{\sigma}_t^{\delta}(\theta)}\right)\frac{\partial^2 \sigma_t^{\delta}(\theta)}{\partial\theta_i \partial\theta_j} \left[\frac{R_t^2}{\sigma_t^{2}(\theta)}-1\right] }_{ =:G_{5,t}^{ij}(\theta)} \\
 && + \frac{2}{n\delta} \sum\limits_{t=1}^{n} \underbrace{\left(\frac{\partial^2 \sigma_t^{\delta}(\theta)}{\partial\theta_i\partial\theta_j} - \frac{\partial^2 \tilde{\sigma}_t^{\delta}(\theta)}{\partial\theta_i\partial\theta_j}\right) \frac{1}{\tilde{\sigma}_t^{\delta}(\theta)} \left[\frac{R_t^2}{\sigma_t^{2}(\theta)}-1\right]}_{=:G_{6,t}^{ij}(\theta)}  \\
 && + \frac{2}{n\delta} \sum\limits_{t=1}^{n} \underbrace{  \left(\frac{R_t^2}{\sigma_t^{2}(\theta)}- \frac{R_t^2}{\tilde{\sigma}_t^{2}(\theta)}\right) \frac{1}{\tilde{\sigma}_t^{\delta}(\theta)}\frac{\partial^2 \tilde{\sigma}_t^{\delta}(\theta)} {\partial\theta_i\partial\theta_j} }_{=:G_{7,t}^{ij}(\theta)}   \ .
\end{eqnarray*}
To prove the desired result we will among others use the following properties $(i): \ \sup\limits_{\theta\in V(\theta_0)}(1/\sigma_t^{\delta}) \leq 1/\underline{\omega}$, $(ii):$ Cauchy Schwartz inequality, $(iii):$ Jensen inequality, $(iv)$: $\partial^2 \tilde{\sigma}_t^{\delta}(\theta)/ \partial\theta_i\partial\theta_j$ is $\alpha/4$ stationary, $(v)$: $\partial \tilde{\sigma}_t^{\delta}(\theta)/ \partial\theta_i$ is $\alpha/2$ stationary.  Let us study $G_{k,t}^{ij}(\theta), \ k=1,...,7$ separately\footnote{we will not study in detail as previously. A detailed proof can be done similarly as in the first part of this lemma.}.

For $G_{1,t}^{ij}(\theta)$ we use $(i)$, $(ii)$, $(iii)$, $(v)$, \eqref{eq:difference between observed and unobserved volatilities} to obtain $E_{\theta_0}\sup\limits_{\theta\in V(\theta_0)}\left|G_{1,t}^{ij}(\theta)\right|^{\alpha/8}  \leq  \mathcal{O}\left(\left[\rho_1^{1/8}\right]^t\right)$.

For $G_{2,t}^{ij}(\theta)$ and $G_{3,t}^{ij}(\theta)$ we use $(i)$, $(ii)$, $(iii)$, $(v)$, \eqref{eq: diff observ and unobserv first der volatility} to obtain $E_{\theta_0}\sup\limits_{\theta\in V(\theta_0)}\left|G_{2,t}^{ij}(\theta)\right|^{\alpha/8}  \leq  \mathcal{O}\left(\left[\rho_3^{1/4}\right]^t\right)$ and $E_{\theta_0}\sup\limits_{\theta\in V(\theta_0)}\left|G_{3,t}^{ij}(\theta)\right|^{\alpha/8}  \leq  \mathcal{O}\left(\left[\rho_3^{1/4}\right]^t\right)$.

For $G_{4,t}^{ij}(\theta)$ we use $(i)$, $(ii)$, $(iii)$ $(v)$ and \eqref{eq: diff observ and unobserv second der volatility} to obtain $E_{\theta_0}\sup\limits_{\theta\in V(\theta_0)}\left|G_{4,t}^{ij}(\theta)\right|^{v/4}  \leq  \mathcal{O}\left(\left[\rho_5^{1/4}\right]^t\right)$.

For $G_{5,t}^{ij}(\theta)$ we use $(i)$, $(ii)$, $(iii)$, $(iv)$, \eqref{eq:difference between observed and unobserved volatilities} and we get $E_{\theta_0}\sup\limits_{\theta\in V(\theta_0)}\left|G_{5,t}^{ij}(\theta)\right|^{\alpha/8} \leq  \mathcal{O}\left(\left[\rho_1^{1/8}\right]^t\right)$.

For $G_{6,t}^{ij}(\theta)$ we use $(i)$, $(ii)$, $(iii)$, \eqref{eq: diff observ and unobserv second der volatility} and we get
$E_{\theta_0}\sup\limits_{\theta\in V(\theta_0)}\left|G_{6,t}^{ij}(\theta)\right|^{\alpha/8} \leq  \mathcal{O}\left(\left[\rho_4^{1/2}\right]^t\right)$.

For $G_{7,t}^{ij}(\theta)$ we use $(i)$, $(ii)$, $(iii)$, $(iv)$ and \eqref{eq: diff observ and unobserv second der volatility} and we get
$E_{\theta_0}\sup\limits_{\theta\in V(\theta_0)}\left|G_{7,t}^{ij}(\theta)\right|^{v/4}  \leq  \mathcal{O}\left(\left[\rho_5^{1/4}\right]^t\right)$.

Using these calculations together with the Jensen inequality and the fact that $\min\{v/4, \alpha/8\}=v/4<1$ we then obtain
\begin{eqnarray*}
E_{\theta_0}\left(n^{v/4} \sup\limits_{\theta\in V(\theta_0)} \left\|H_n(\theta)-\tilde{H}_n(\theta)\right\|^{v/4}\right) &=&
E_{\theta_0}\sup\limits_{\theta\in V(\theta_0)} \left( n \left\|H_n(\theta)-\tilde{H}_n(\theta)\right\|\right)^{v/4}\\
&=& E_{\theta_0}\sup\limits_{\theta\in V(\theta_0)} \left( n \sum\limits_{i=1}^{m}\sum\limits_{j=1}^{m} \left|\left(H_n(\theta)-\tilde{H}_n(\theta)\right)_{ij}\right|\right)^{v/4}\\
&\leq& E_{\theta_0} \sup\limits_{\theta\in V(\theta_0)} n^{v/4} \sum\limits_{i=1}^{m}\sum\limits_{j=1}^{m}   \left|\left(H_n(\theta)-\tilde{H}_n(\theta)\right)_{ij}\right|^{v/4} \\
&=& \sum\limits_{i=1}^{m}\sum\limits_{j=1}^{m} E_{\theta_0} \sup\limits_{\theta\in V(\theta_0)} n^{v/4} \left|\left(H_n(\theta)-\tilde{H}_n(\theta)\right)_{ij}\right|^{v/4} \\
&\leq&  \sum\limits_{i=1}^{m}\sum\limits_{j=1}^{m} E_{\theta_0} \sup\limits_{\theta\in V(\theta_0)} \left(\frac{2}{\delta}\right)^{v/4} \left|\sum\limits_{k=1}^{7} \sum\limits_{t=1}^{n} G_{k,t}^{ij} \right|^{v/4} \\
&\leq&  \left(\frac{2}{\delta}\right)^{v/4} \sum\limits_{i=1}^{m}\sum\limits_{j=1}^{m} \sum\limits_{k=1}^{7} \sum\limits_{t=1}^{n} E_{\theta_0} \sup\limits_{\theta\in V(\theta_0)}  \left| G_{k,t}^{ij} \right|^{v/4} \\
&\leq&  \left(\frac{2}{\delta}\right)^{v/4} \sum\limits_{i=1}^{m}\sum\limits_{j=1}^{m} \sum\limits_{k=1}^{7} \underbrace{ \sum\limits_{t=1}^{\infty} E_{\theta_0} \sup\limits_{\theta\in V(\theta_0)}  \left| G_{k,t}^{ij} \right|^{v/4}}_{<\infty \textrm{ from above calculations}} \\
&=:& F_2  \ ,
\end{eqnarray*}
where $F_2<\infty$ and independent of $n$.

From this we can obtain the required convergence result by using the Markov inequality:
\begin{eqnarray*}
\mathbb{P}\left(\sup\limits_{\theta\in V(\theta_0)} \left\|H_n(\theta)-\tilde{H}_n(\theta)\right\| > \varepsilon \right) &=& \mathbb{P}\left(n^{v/4}\sup\limits_{\theta\in V(\theta_0)} \left\|H_n(\theta)-\tilde{H}_n(\theta)\right\|^{v/4} > n^{v/4}\varepsilon^{v/4} \right) \\
&\leq & \frac{E_{\theta_0}\left(n^{v/4} \sup\limits_{\theta\in V(\theta_0)} \left\|H_n(\theta)-\tilde{H}_n(\theta)\right\|^{v/4}\right)}{n^{v/4}\varepsilon^{v/4}} \\
&=& \frac{F_2}{n^{v/4}\varepsilon^{v/4}} \ \longrightarrow \ 0 \  (as \  n \rightarrow \infty) \ \textrm{for every } \varepsilon>0.
\end{eqnarray*}
This proves the second equation of the lemma. $\hspace{7cm} \square$
\newline

 
\textbf{\textit{Proof of the proposition \ref{prop: ass norm of QMLE}}} \newline
As shown previously, it suffices to prove the equations \eqref{eq: first main equation ass norm} and \eqref{eq: second main equation ass norm}.

From the lemma \ref{lemma: second lemma for ass norm}, the process $\partial\ell_t(\theta_0)/\partial\theta$ is stationary and ergodic. Moreover it is a martingale difference sequence since $E_{\theta_0}(\partial\ell_t(\theta_0)/\partial\theta| \mathfrak{F}_{t-1})=0$ as it can be seen in the proof of lemma \ref{lemma: first lemma for ass norm}. Furthermore from this lemma \ref{lemma: first lemma for ass norm}, we have that this process has finite variance. Hence using the invariance principle of stationary martingale difference, we obtain
\begin{equation*}
\frac{1}{\sqrt{n}}\sum\limits_{t=1}^{n}{\frac{\partial\ell_t(\theta_0)}{\partial\theta}} \ \longrightarrow \ \mathcal{N}\left(0, \frac{4}{\delta^2}(\kappa-1)A\right)  \ \textrm{(in distribution),}
\end{equation*}
or equivalently,
\begin{equation*}
\sqrt{n}S_n(\theta_0) \ \longrightarrow \ \mathcal{N}\left(0, \frac{4}{\delta^2}(\kappa-1)A\right)  \ \textrm{(in distribution).}
\end{equation*}
Since the lemma \ref{lemma: fourth lemma for ass norm diff ell and tilde ell} says that the difference between the observed $(\tilde{S}_n)$ and unobserved $(S_n)$ volatilities converges to zero in probability, we can use it together with the above equation to obtain through the Slutsky lemma that
\begin{equation*}
\sqrt{n}\tilde{S}_n(\theta_0) = \sqrt{n}S_n(\theta_0) + \sqrt{n}\left(\tilde{S}_n(\theta_0)-S_n(\theta_0)\right) \ \longrightarrow \ \mathcal{N}\left(0, \frac{4}{\delta^2}(\kappa-1)A\right) \  \ \textrm{(in distribution).}
\end{equation*}
This proves the equation \eqref{eq: first main equation ass norm}. 

To prove the equation \eqref{eq: second main equation ass norm} we used the fact that $\partial^2\ell_t/\partial\theta \partial\theta^{'}$ is stationary and ergodic (lemma \ref{lemma: second lemma for ass norm}) and that $E_{\theta_0}\left|\partial^2\ell_t/\partial\theta \partial\theta^{'}\right|<\infty$ to conclude -using the weak law of large numbers- that for all $\theta \in V^{'}(\theta_0)$,
\begin{equation}
H_n(\theta) \  \longrightarrow \ H(\theta) \ \textrm{(in prob.)} \  \ as \ n\rightarrow\infty  \ ,
\label{pointwise convergence oh Hn}
\end{equation}
where $H(\theta)=E_{\theta_0} \left(\partial^2\ell_t(\theta)/\partial\theta \partial\theta^{'}\right).$

Furthermore the third derivative of $\ell_t$ is bounded in $V^{'}(\theta_0)$. Hence the above pointwise convergence is uniform over $V^{'}(\theta_0)$, i.e.,
\begin{equation*}
\sup\limits_{\theta\in V^{'}(\theta_0)} \left\|H_n(\theta)-H(\theta)\right\| \ \longrightarrow \ 0 \ \textrm{(in prob.)} \  \ as \ n\rightarrow\infty.
\end{equation*}

Combining this with the equation \eqref{eq:diff Hn and tilde Hn} of the lemma \ref{lemma: fourth lemma for ass norm diff ell and tilde ell}, we get
\begin{equation}
\sup\limits_{\theta\in V^{'}(\theta_0) \cap V(\theta_0)} \left\|\tilde{H}_n(\theta)-H(\theta)\right\| \ \longrightarrow \ 0 \ \textrm{(in prob.)} \  \ as \ n\rightarrow\infty.
\label{uniform convergence oh Hn2}
\end{equation}

From this uniform convergence of $H_n$ towards $H$ over $ V^{'}(\theta_0) \cap V(\theta_0)$ and the fact that $\tilde{H}_n$ is continuous in $\theta$, we obtain that $H$ is also continuous over $ V^{'}(\theta_0) \cap V(\theta_0)$. Using now the continuity of $H$ and the fact that that $\bar{\theta}_n$ converges towards $\theta_0$ as $\hat{\theta}_n$ converges towards $\theta_0$, together with the above pointwise convergence \eqref{pointwise convergence oh Hn}, we get
\begin{eqnarray*}
\left\|H_n(\bar{\theta}_n)-H(\theta_0)\right\| &\leq& \left\|H_n(\bar{\theta}_n)-H(\bar{\theta}_n)\right\| + \left\|H(\bar{\theta}_n)-H(\theta_0)\right\|  \longrightarrow \ 0 \ \textrm{(in prob.)} \  \ as \ n\rightarrow\infty \ .
\end{eqnarray*}
Furthermore
\begin{eqnarray*}
H_{i,j}(\theta_0) &=& E_{\theta_0} \left(\frac{\partial^2 \ell_t(\theta_0)}{\partial\theta_i\partial\theta_j}\right) \\
&=& E_{\theta_0} \left( \frac{2}{\delta} \frac{1}{\sigma_t^{\delta}(\theta_0)} \left[-\frac{1}{\sigma_t^{\delta}(\theta_0)} \frac{\partial \sigma_t^{\delta}(\theta_0)}{\partial\theta_i} \frac{\partial \sigma_t^{\delta}(\theta_0)}{\partial\theta_j} \left(\frac{2+\delta}{\delta} \frac{R_t^2}{\sigma_t^{2}(\theta_0)}-1 \right)  +  \frac{\partial^2 \sigma_t^{\delta}(\theta_0)}{\partial\theta_i \partial\theta_j} \left( \frac{R_t^2}{\sigma_t^{2}(\theta_0)}-1 \right)  \right]\right) \\
&=& E_{\theta_0} \left( \frac{-2}{\delta} \frac{1}{\sigma_t^{2\delta}(\theta_0)} \frac{\partial \sigma_t^{\delta}(\theta_0)}{\partial\theta_i} \frac{\partial \sigma_t^{\delta}(\theta_0)}{\partial\theta_j}\right) \underbrace{E_{\theta_0} \left(\frac{2+\delta}{\delta} \varepsilon_t^2-1 \right)}_{=2/\delta}  +  E_{\theta_0} \frac{\partial^2 \sigma_t^{\delta}(\theta_0)}{\partial\theta_i \partial\theta_j} \underbrace{E_{\theta_0} \left(\varepsilon_t^2-1 \right)}_{=0} \\
&=& \frac{-4}{\delta^2} A_{i,j}.
\end{eqnarray*}
This proves the proposition.   $\hspace{7cm} \square$
\section*{Acknowledgment}
 The authors thank Peter Ruckdeschel for his ideas and comments. The authors also thank the financial support of the Fraunhofer Institut f\"ur Techno- und Wirtschaftsmathematik.

\newpage
\bibliographystyle{econ_letters}
\bibliography{myrefs}

\end{document}